\documentclass[onefignum,onetabnum,onethmnum]{siamonline220329}
\usepackage{amsfonts,amssymb}
\usepackage{graphicx}
\usepackage{epstopdf}
\usepackage{enumitem}
\usepackage{cite}
\setlist[enumerate]{leftmargin=.5in}
\setlist[itemize]{leftmargin=.5in}

\crefname{hypothesis}{Hypothesis}{Hypotheses}
\newsiamthm{claim}{Claim}
\renewtheorem{proposition}{Proposition}
\renewtheorem{lemma}{Lemma}
\newtheorem{remark}{Remark}
\newtheorem{example}{Example}

\usepackage{color}
\usepackage{algpseudocode}

\usepackage{epstopdf}
\usepackage{dsfont}

\usepackage{color}

\newcommand{\one}{{\mathbf{1}}}

\newcommand{\tran}{^{\top}}

\newcommand{\diag}{\mbox {\rm diag}}

\newcommand{\beq}{\begin{equation}}
\newcommand{\eeq}{\end{equation}}
\newcommand{\bea}{\begin{eqnarray}}
\newcommand{\eea}{\end{eqnarray}}
\newcommand{\beas}{\begin{eqnarray*}}
\newcommand{\eeas}{\end{eqnarray*}}
\newcommand{\ba}{\begin{array}}
\newcommand{\ea}{\end{array}}
\newcommand{\bit}{\begin{itemize}}
\newcommand{\eit}{\end{itemize}}
\newcommand{\ben}{\begin{enumerate}}
\newcommand{\een}{\end{enumerate}}

\newcommand{\ped}[1]{ _{ {\mathrm{#1} } }}

\newcommand{\ap}[1]{ ^{ {\mathrm{#1} } }}

\newcommand{\Real}[1]{ { {\mathbb R}^{#1} } }

\newcommand{\calD}{{\cal D}}
\newcommand{\calG}{{\cal G}}

\newcommand{\calI}{{\cal I}}
\newcommand{\calE}{{\cal E}}

\newcommand{\calS}{{\cal S}}

\newcommand{\calV}{{\cal V}}

\definecolor{lgray}{gray}{0.7}

\def\ve{\varepsilon}

\headers{Optimal Clearing Payments in a Financial Contagion Model}{G. Calafiore, G. Fracastoro, A.V. Proskurnikov}

\title{Optimal Clearing Payments in a Financial Contagion Model\thanks{%First submitted to the editors on June 13, 2022. 
A preliminary version of this work~\cite{Calafiore2021a} was presented on the 60th IEEE Conference on Decision and Control, December 13-15, 2021.
\funding{This study was carried out within the project 2022K8EZBW ``Higher-order interactions in social dynamics with application to monetary networks'', funded by European Union – Next Generation EU  within the PRIN 2022 program (D.D. 104 - 02/02/2022 Ministero dell'Universit\'a e della Ricerca). This manuscript reflects only the authors’ views and opinions, and the Ministry cannot be considered responsible for them.
}
}}

\author{Giuseppe C. Calafiore\thanks{DET - Politecnico di Torino, Italy, \email{giuseppe.calafiore@polito.it}, IEIIT-CNR Italy, and VinUniversity, Hanoi, Vietnam, \email{giuseppe.c@vinuni.edu.vn}}
\and Giulia Fracastoro\thanks{DET - Politecnico di Torino, \email{giulia.fracastoro@polito.it}}
\and Anton V. Proskurnikov\thanks{DET - Politecnico di Torino, \email{anton.p.1982@ieee.org}}
}

\usepackage{amsopn}

\begin{document}
\maketitle

\begin{abstract}
Financial networks are characterized by complex structures of mutual obligations. These obligations are fulfilled entirely or in part (when defaults occur) via a mechanism called {\em clearing}, which  determines a set of payments that settle the claims by respecting rules such as limited liability, absolute priority, and proportionality (pro-rated payments).
In the presence of shocks on the financial system, however, the clearing mechanism may lead to cascaded defaults and eventually to financial disaster.
In this paper, we first study the clearing model under pro-rated payments of Eisenberg and Noe, and we derive novel  necessary and sufficient conditions for the uniqueness of the clearing payments,
valid for an arbitrary topology of the financial network. Next, we observe that
the proportionality rule is a factor that potentially concurs to the cascaded defaults effect, and that
the aggregated systemic loss can be reduced if this rule is lifted. We thus  shift the focus from the individual interest to the overall systemic interest to contain the adverse effects of cascaded failures, and we show that pro-rate-free clearing payments  can be computed  uniquely by solving suitable convex optimization problems.
\end{abstract}

%\MSCCLASS{}
%\ORMSCLASS{Primary: ; secondary: }
%\HISTORY{}
\section{Introduction}

The financial industry is participated by organizations that are
linked to each other by means of  an intricate structure of mutual obligations.
The behavior of such financial interconnected system has been extensively studied over the past years, see for instance
\cite{Gale2007,Battiston2010}. The interconnection among financial institutions creates potential channels of contagion, whereby a failure (financial default) of a single entity in  the system can result in a threat to the stability of the entire global financial system.
Recent examples of such a behavior include the collapse of Lehman Brothers, recognized as one of the reasons of the global financial crisis in 2008, and the government bailout of the giant insurance company AIG,~\cite{BlundellWignall2010}.
Much effort has  been invested in understanding the effects of {\em systemic risk}, that is,  how stresses, such as bankrupts and failures, to one part of the system can spread to others, and eventually lead to avalanche breakdowns, see, e.g.,~\cite{Eisenberg2001,Elliott2014,Glasserman2016,Haldane2011}.

An important line of research pursued in systemic risk theory focuses on the development of realistic models of \emph{clearing} procedures between financial institutions.
Clearing is essentially
a set of rules under which the participants to a financial network agree to settle payments, when these payments cannot
meet the original liabilities due to defaults,~\cite{Kabanov2018}.
The seminal work in~\cite{Eisenberg2001} introduced a simple mathematical model of clearing in a financial network, in which financial institutions have two types of assets: the external assets (e.g., incoming cash flows) and the internal assets (e.g.,  funds that banks lend to one another).
The model in~\cite{Eisenberg2001} assumes that the obligations of all entities within the financial system are paid simultaneously and are determined by  three fundamental rules:
1) limited liability, that is, the total payment of each node can not exceed its available cash flow;
2) the priority of the debt claims, that is,   stockholders  receive no value until the node is able to completely pay off all of its outstanding liabilities;
3) the proportionality, or pro-rata rule, that is, all debts have equal priority, so that all claimant institutions are paid proportionally to their nominal claims.
Under these assumptions, the matrix of mutual interbank payments is uniquely determined by the so-called \emph{clearing vector}, which is found as the fixed point of a nonlinear equation. This vector always exists~\cite{Eisenberg2001}, whereas its uniqueness has been proved under certain regularity assumptions, see~\cite{Eisenberg2001,Glasserman2016,Kabanov2018,Rogers2013}; these uniqueness conditions are however only \emph{sufficient}, but not necessary.

The basic model offered in~\cite{Eisenberg2001} has been later extended in various directions, incorporating non-trivial features of real-world financial networks, see, e.g., the recent survey in~\cite{Jackson2021}. The models presented in~\cite{Suzuki2002,Elsinger2009}, for instance, take into account cross-holdings and seniorities of liabilities, and the works of~\cite{Cifuentes2005,Shin2008} introduce the concept of liquidity risk. Other works considered also illiquid assets \cite{Amini2016}, cross-ownership of equities and liabilities~\cite{Fischer2014}, decentralized clearing processes~\cite{Csoka2018}, and multiple maturity dates~\cite{Kusnetsov2019}.

The contribution of the present work is twofold. First, we give a full solution to the problem of the uniqueness of the clearing vector in the proportional payments case.
From the financial perspective, uniqueness is important since it guarantees that no ambiguity exists in the payments, so that each entity must abide to one and only one clearing payment, with no possible controversy.
From the computational viewpoint, uniqueness ensures that different methods for finding a clearing vector,
such, e.g., the fictitious default algorithm~\cite{Glasserman2016} and more advanced methods from~\cite{Kusnetsov2019}, return the same answer.
The first sufficient graph-theoretical condition for the uniqueness of the clearing payments was obtained in~\cite{Eisenberg2001}
(see also the works~\cite{ElBitar2017,Kabanov2018}, giving a simpler and more elegant proof): the clearing vector is unique if the financial network is \emph{regular}, which means that every bank either has an outside asset, or has a (direct or indirect) creditor with outside assets.
Another sufficient condition for uniqueness  is formulated in~\cite{Glasserman2016}: the clearing vector is unique if each node of the network has a chain of liabilities to the external sector. Both conditions are only sufficient yet not necessary.
%Under the assumption of a network's \emph{strong connectivity}, a necessary and sufficient condition for the uniqueness was found in~\cite{REN2016779}.
To the best of our knowledge, the only necessary and sufficient condition for the clearing vector's uniqueness applicable to an arbitrary financial network available in the literature is the very general result in~\cite{Massai2022}, which examines the uniqueness of equilibria in a dynamical flow network with saturations.
This criterion, primarily motivated by more general models of systemic risk proposed in~\cite{Elsinger2006}, however,
appears to be of limited practical  use in the classical Eisenberg-Noe model, since it requires computation of some parameters that depend on the payment matrix (left and right Perron eigenvectors for its irreducible blocks). Unlike the procedure in~\cite{Massai2022}, our proposed method allows to test the uniqueness of the clearing vector without knowledge of the payment matrix; only the graph of liability relations matters. Similar to~\cite{Massai2022}, we are also able to find the whole set of clearing vectors. As a byproduct of the developed theory, we also derive some new properties of the maximal (or dominant) clearing vector that are of independent interest.

The second key aspect addressed in this paper is the lifting of the proportionality (pro-rata) rule.
Such division rule is the most common in the literature on systemic risk and it was introduced since the seminal work of~\cite{Eisenberg2001}. The proportional rule, however, may not be realistic in some situations, for instance when liabilities have different seniorities, see, e.g.,~\cite{Kaminski2000}. Other division rules have thus been also considered in the literature; examples include %priority rules~\cite{Moulin2000,Chatterjee2015}, 
constrained equal awards or constrained equal losses rules~\cite{Thomson2013}, welfare-maximizing rules~\cite{Gallice2019}, and general division rules \cite{Csoka2018,Barratt2020}. In~\cite{Csoka2018}, the authors propose a decentralized clearing process in a discrete setup which assumes integer payments. Instead,~\cite{Barratt2020} considers a multi-period liability clearing problem. Differently from the model considered in this work, \cite{Barratt2020} also makes the assumption that entities cannot pay other entities more than the cash they have on hand.
Motivated by these works, in Section~\ref{sec:no_prorata} a clearing scheme under unconstrained (i.e., non necessarily proportional) payments is analyzed. We show that an optimal clearing matrix  can be computed by solving a linear program.
Plain relaxation of the proportionality constraint leads as a downside to the  loss of the clearing matrix uniqueness.
Also, the feasible payment matrices do not constitute a complete lattice, in particular, even for very simple financial networks with $n=3$ nodes the minimal and maximal payment matrix may fail to exist.
However, we recover uniqueness by means of a simple two-stage optimization scheme whose result is the unique clearing payments matrix that {\em (a)} achieves the best possible deviation loss from the nominal liabilities, and {\em (b)} has the minimum ``size'' among all matrices achieving the optimal loss.
The positive systemic effects  obtained by releasing the pro-rata rule in favor of an aggregated performance approach are illustrated by
means of a schematic example as well as via extensive numerical tests  using a synthetic random network, which is similar to a ``testbench'' model from~\cite{Nier2007}.

The practical implementation  of an optimal unconstrained clearing payments scheme, however, currently faces several obstacles. Some of these obstacles are actually shared also by the pro-rated scheme and by other mathematical default schemes, and are due to unmodelled non-idealities, contract renegotiations, credit freeze, local laws and delays in their application. Another issue is that any clearing scheme should be contractualized ex-ante and all players should adhere to it.  More relevantly, an optimal aggregated unconstrained clearing approach presupposes the existence of a central authority that has full knowledge of the inter-bank liability matrix. While this assumption may not be valid in the current situation, we believe that the technology exists for making the financial network more transparent, and
that evidence of the systemic and social advantages of an aggregated approach to clearing may foster the progress in this area.

The remainder of this paper is organized as follows. Section~\ref{sec:notation} defines the notation.
Section~\ref{sec:fin_nets} introduces the Eisenberg-Noe model and related concepts.
Section~\ref{sec:prorata} presents a novel necessary and sufficient condition for the clearing vector's uniqueness, in the situation where the pro-rata constraint is adopted.
 In Section~\ref{sec:no_prorata}, we study the set of clearing matrices who do not necessarily satisfy the pro-rata rule. We demonstrate that one such matrix can always be found by solving a convex optimization problem aimed at minimizing a system-level loss.
 In general, however,  the optimal clearing matrix is non-unique. We next consider specifically a system-level (i.e., aggregated) loss given by the sum of all deviations of the actual payments from the nominal liabilities. In this setting, we provide a two-stage procedure based on convex programming for finding the unique clearing payments matrix which achieves the best possible system loss and which has at the same time the smallest possible Euclidean size.
Section~\ref{sec:experiments} presents a schematic example that shows
how the proposed approach may lead to effective isolation of shocks and containment of the default contagion.
This section also contains extensive  simulation tests on random networks, with a comparison
of the system losses and the number of default nodes in the cases where the proportionality rule is adopted and  where it is discarded. Section~\ref{sec:conclusion} concludes the paper.

\section{Preliminaries and notation}
\label{sec:notation}
Given a finite set $\calV$, the symbol $|\calV|$ stands for its cardinality.
For two families of real numbers $(a_{\xi})_{\xi\in\Xi},(b_{\xi})_{\xi\in\Xi}$, the symbol $a\leq b$ ($b$ dominates $a$, or $a$ is dominated by $b$) denotes the element-wise relation $a_{\xi}\leq b_{\xi}$, $\forall\xi\in\Xi$. We write
$a\lneq b$ if $a\leq b$ and $a\ne b$. The operations $\min$ and $\max$ are defined element-wise, e.g., $\min(a,b)\doteq(\min(a_{\xi},b_{\xi}))_{\xi\in\Xi}$. These notation symbols apply to vectors (usually, $\Xi=\{1,\ldots,n\}$) and matrices (usually, $\Xi=\{1,\ldots,n\}\times\{1,\ldots,n\}$).

Every nonnegative square matrix $A=(a_{ij})_{i,j\in\calV}$ corresponds to a (weighted directed) graph $\calG[A]=(\calV,\calE[A],A)$ whose nodes are indexed by $\calV$ and whose set of arcs is defined as $\calE[A]=\{(i,j):a_{ij}>0\}$. The value $a_{ij}$ can be interpreted as a weight of arc $(i,j)$, which is also denoted as $i\rightarrow j$. A sequence of arcs $i_0\rightarrow i_1\rightarrow\ldots\rightarrow i_{s-1}\rightarrow i_s$ constitutes a \emph{walk} between nodes $i_0$ and $i_s$ in graph $\calG[A]$. The set of nodes $J\subseteq\calV$  is \emph{reachable} from node $i$ if $i\in J$ or a walk from $i$ to some element $j\in J$ exists; $J$ is called \emph{globally reachable} in the graph if it is reachable from every node $i\not\in J$.

A graph is strongly connected (strong) if every two nodes $i,j$ are mutually reachable. A graph that is not strong has several strongly connected components (for brevity, we call them simply \emph{components}). A component is said to be non-trivial if it contains more than one node. A component is said to be a \emph{sink} component if no arc leaves it. %and a \emph{source} component if no arc enters it.

A nonnegative square matrix $A\in\Real{\calI\times\calI}$ is said to be \emph{stochastic} if all its rows sum to $1$: $\sum_{j\in\calI} a_{ij}=1$, $\forall i\in\calI$ and \emph{substochastic} if $\sum_{j\in\calI} a_{ij}\leq 1$, $\forall i\in\calI$. Introducing the vector of ones $\one\in\Real{\calI}$, matrix $A\geq 0$ is stochastic if $A\one=\one$ and substochastic if $A\one\leq\one$.

\section{Financial Networks}
\label{sec:fin_nets}

We henceforth  use the notation introduced in~\cite{Glasserman2016}, except for a few minor changes.
A \emph{financial network} is represented as a weighted digraph $\cal G=(\calV,\calE,\bar P)$ whose
nodes stand for financial institutions (banks, funds, insurance companies, etc.) and whose weighted adjacency matrix $\bar P=(\bar p_{ij})$ represents the mutual liabilities among the institutions. Namely, entry $\bar p_{ij}\geq 0$ means
that node $i$ has an obligation to pay $\bar p_{ij}$ currency units to
 node $j$ at the end of the current time period, and an arc $(i,j)\in\calE$ from node $i$ to node $j$ exists if and only if $\bar p_{ij}>0$.
By definition, $\bar p_{ii}=0\,\forall i$, so the graph contains no self-arcs.

Along with mutual liabilities, the banks have \emph{outside assets}. The outside \emph{asset} $\bar c_i\geq 0$ is the total payment due from non-financial entities (the external sector) to node $i$; these numbers constitute vector $\bar c=(\bar c_i)_{i\in\calV}$.
%Similarly, one can consider the  liability of node $i$
%towards the external sector
%as $\bar b_i\geq 0$. Often, the external liabilities are formally replaced by liabilities  to an additional ``fictitious'' node, representing the external sector~\cite{Eisenberg2001}.
%Adding this ``virtual'' node to $\calV$, we shall henceforth assume without loss of generality that $\bar b=0$.

For each node $i$, we define the nominal cash in-flow and out-flow (standing for the \emph{asset} and \emph{liability} sides of the balance sheet):
\begin{equation}\label{eq.in-out-nominal}
\bar\phi_i\ap{in}\doteq  \bar c_i + \sum_{k\neq i} \bar p_{ki},\quad
\bar p_i \doteq  \bar\phi_i\ap{out}\doteq   \sum_{k\neq i} \bar p_{ik}.
\end{equation}
\begin{remark}\rm
Notice that the network may contain one or several \emph{sink nodes}, that is nodes without outgoing arcs ($\bar p_i=0$). These nodes represent banks with assets but without liabilities. One sink node is  often introduced (see, e.g.,~\cite{Eisenberg2001,Glasserman2016}) for accommodating liabilities to non-financial institutions: it represents a  fictitious  financial entity with  no liabilities and whose assets are the liabilities of the remaining banks to the external sector.
\end{remark}

%The nodes with $\bar p_i=0$ have no outgoing arcs and, according to the  graph-theoretical terminology, they are called \emph{sinks}.
%As mentioned before,  one such node can be fictitiously defined for the purpose of collecting the debts to the external sector. In general,   however, other sinks %may exist.

In regular operations it will hold that $\bar\phi_i\ap{in}\geq \bar\phi_i\ap{out}$, meaning that each bank is able to pay its debts at the end of the  period. The risk of financial contagion arises in the situation when
a financial shock hits some nodes, meaning that the outside assets drop to smaller-than-expected values
$c_i\in [0,\bar c_i)$. In this situation, it may happen that
\[
c_i + \sum_{k\neq i} \bar p_{ki}<\bar p_i.
\]
In this case, node $i$ becomes unable to fully meet its payment obligations, and then  \emph{defaults}.
When in default, a node pays out according to its capacity, thus reducing
the amounts paid to the adjacent nodes, which in turn, for this reason, may also default and reduce their payments to other nodes, and so on in a cascaded fashion.
As a result of default, the \emph{actual} payment $p_{ij}\in [0,\bar p_{ij}]$ from node $i$ to node $j$, in general, may be less than the nominal due payment $\bar p_{ij}$. A natural question arises: which matrices of actual payments $P=(p_{ij})\leq\bar P$ may be considered as ``fair''
in the case of default?
We shall see that the pro-rata rule is a commonly accepted rule for allocating payments in the case of default, but we shall also explore an alternative approach that aims at minimizing the aggregated loss over the financial system in Section~\ref{sec:no_prorata}.
Denote the vectors of actual in-flows and out-flows by
\begin{equation}\label{eq:in-out-flows}
\phi\ap{in} \doteq c + P\tran \one,\quad
 p\doteq\phi\ap{out}\doteq P \one.
\end{equation}
The conditions to which the  payments matrix $P=P(c,\bar P)$  is subject to are as follows,~\cite{Eisenberg2001}:
\begin{enumerate}[label=(\roman*)]
\item \textbf{(limited liability)} The total payment of each node does not exceed its in-flow, that is, $\phi\ap{in}\geq\phi\ap{out}$;
\item \textbf{(absolute priority of debt claims)} Either node $i$ pays its obligations in full ($p_i=\bar p_i$), or it pays all its value to the creditors ($p_i=\phi_i\ap{in}$).
\end{enumerate}
Recalling that $P\leq \bar P$ and $p=P\one\leq\bar p$, conditions (i) and (ii) are reformulated compactly as
\beq\label{eq:clearing-1}
P\one=\min(\bar p,c + P\tran \one).
\eeq
\begin{definition}\label{def.matrix}
A matrix $P$ is called a \emph{clearing matrix} (or matrix of clearing payments) corresponding to the vector of outside assets $c$, if $0\leq P\leq\bar P$, and~\eqref{eq:clearing-1} holds.
\end{definition}

\vspace{.2cm}
Notice that~\eqref{eq:clearing-1} is a system of $n\doteq|\calV|$ nonlinear equations in $n^2$ variables $p_{ij}$. Hence, one cannot expect to find a unique solution, in general. To obtain uniqueness of the solution (in the generic situation), a third requirement is typically introduced, see, e.g.,~\cite{Eisenberg2001}, known as the \emph{proportionality} or \emph{pro-rata} rule, which expresses the
requirement that all debts have equal priority and must be paid in proportion to the initial claims.
The imposition of this rule reduces the number of variables to $n=|\calV|$. It is known that under the pro-rata rule a clearing vector always exists;
also, one such vector can be found by solving a convex optimization problem with $n$ variables, applying a standard fixed-point iteration or a more advanced ``fictitious default algorithm''~\cite{Eisenberg2001,Glasserman2016,Kusnetsov2019}.
In Section~\ref{sec:prorata} we
provide a necessary and sufficient criterion for the uniqueness of the clearing vector in the pro-rata case, and we also offer an algorithm for finding the set of all clearing vectors.
The pro-rata rule reflects an underlying criterion of { local fairness} among neighboring nodes, and it is a convention enforced in many contracts. In Section~\ref{sec:no_prorata}
we discuss  the case where the pro-rata  rule is lifted
and substituted by a system-level  aggregated loss minimization criterion.

\section{Pro-rata rule and clearing vectors}\label{sec:prorata}

One standard approach to determine the clearing payments is based on imposing an additional restriction on
the payments $p_{ij}$,
 stating that the payments  of node $i$ to the claimants should be proportional to the nominal liabilities $\bar p_{ij}$.
 It is convenient to introduce the matrix of normalized, or \emph{relative} liabilities
\begin{equation}\label{eq.A}
A=(a_{ij}),\quad
a_{ij}=
\begin{cases}
\frac{\bar p_{ij}}{\bar p_i},& \mbox{if }\bar p_i>0,\\
1, &\mbox{if } \bar p_i=0 \mbox{ and } i=j,\\
0, &\text{otherwise.}
\end{cases}
\end{equation}
By definition, matrix $A$ is \emph{stochastic}, that is, $a_{ij}\geq 0$ and $\sum_j a_{ij}=1$ for all $i$ or, equivalently, $A\one=\one$.
The {pro-rata rule} can then be formulated as $P = \diag(P\one)A$ or, equivalently,
\beq\label{eq:prorata}
p_{ij}=p_ia_{ij},\quad \forall i,j\in\calV.
\eeq
Condition~\eqref{eq:prorata} is known as  equal priority~\cite{Eisenberg2001}, \emph{pro-rata}~\cite{Glasserman2016} or proportionality~\cite{Rogers2013} rule. Under~\eqref{eq:prorata}, it can be proved that
\[
(P^{\top}\one)_i=\sum_{j\in\calV}(P^{\top})_{ij}=\sum_{j\in\calV}p_{ji}=\sum_{j\in\calV}p_ja_{ji}=(A^{\top}p)_i\,\forall i,
\]
which allows us to rewrite~\eqref{eq:clearing-1} in the equivalent vector form
\beq\label{eq:clearing-2}
p=\min(\bar p,c + A^{\top}p).
\eeq

\begin{definition}
%Under the assumption~\eqref{eq:prorata}, the
A vector $p\geq 0$ is said to be a \emph{clearing} vector if it satisfies~\eqref{eq:clearing-2}.
\end{definition}

\vspace{.2cm}
The existence of a clearing vector is usually proved by appealing to the general Knaster-Tarski fixed-point theorem~\cite{Eisenberg2001,Kabanov2018}, applied to the non-decreasing mapping
\[
p\mapsto \min(\bar p,c + A^{\top}p).
\]
This theorem implies that the set of clearing vectors is non-empty and, furthermore, this set constitutes a complete lattice (with respect to the relation $\leq$), therefore, the \emph{minimal} and \emph{maximal} clearing vectors do exist. This monotonicity-based approach is convenient, because it allows to prove the existence of clearing vectors in more complicated models~\cite{Amini2016,Kabanov2018}. 
%The Knaster-Tarski theorem also suggests an iterative procedure to compute some clearing vector, which is also known as the ``fictitious default %algorithm''~\cite{Eisenberg2001}.

At the same time, the Knaster-Tarski theorem does not give a full description of the set of all clearing vectors. One of the important questions that dates back to the original model from~\cite{Eisenberg2001} is whether this set is a singleton, that is, whether the three simple rules (absolute priority, limited liability and pro-rata payments) uniquely determine a clearing vector. Such a uniqueness guarantees that no ambiguity exists
in the payments, so that each entity must abide to one and only one clearing payment, which is
important from both the economical and the computational viewpoints. As it will be shown, in many aspects the maximal (or dominant) clearing vector is the most natural, because it minimizes the overall system loss. The computation of this maximal clearing vector is non-trivial, as usually one has to solve an LP with $n$ variables and $2n$ nonlinear constraints; an alternative method from~\cite{Rogers2013}  finds the maximal clearing vector in no more than $n$ steps,  where at each step one has to solve a non-degenerate system of linear equations of the dimension $O(n)$. If the clearing vector is unique, then it can be computed by a more efficient ``fictitious default algorithm''~\cite{Eisenberg2001,Glasserman2016}.

Note that in degenerate situations the clearing vector may be non-unique. For instance, if $A=I_n$, $\bar p\ne 0$, and $c=0$, every vector $p$ such that $0\leq p\leq \bar p$, obviously, satisfies~\eqref{eq:clearing-2}.
As it will be shown (Theorem~\ref{thm.unique1}),  the existence of such a ``closed'' subgroup of banks independent of the remaining network and external sector is in fact the only reason for non-uniqueness of the clearing vector. %(Theorem~\ref{thm.unique1} below).
Subsection~\ref{subsec.uniq-gen} offers necessary and sufficient conditions for the clearing vector's uniqueness. We also show that, even when these conditions do not hold, still some of the clearing vector's elements are uniquely determined by $A$ and $c$. In such situations, our Theorem~\ref{thm.unique} allows to describe the whole polytope of clearing vectors.

\subsection{The dominant clearing vector -- Extremal properties}

Whereas the existence of a maximal clearing vector is usually proved via the Knaster-Tarski fixed-point theorem~\cite{Eisenberg2001}, we consider an alternative construction, which also clarifies the geometrical meaning of this vector.
Consider the convex polyhedron
\begin{equation}\label{eq.d}
\calD=[0,\bar p]\cap\{p:c+A\tran p\geq p\}.
\end{equation}
The set $\calD$ is non-empty (it contains, e.g., the null vector), and it can thus be represented as  the convex hull of its extreme points (or vertices).
The following lemma (see the appendix Section~\ref{sec:appendix} for its proof)  shows that one of these extreme points is the maximal (with respect to $\leq$ relation) element of $\calD$, being also a clearing vector.

\vspace{.2cm}
\begin{definition}\label{def.decrease}
$F:[0,\bar p]\to\Real{}$ is non-increasing (respectively, decreasing) if $F(p^1)\geq F(p^2)$
(respectively, $F(p^1)>F(p^2)$) whenever $p^1,p^2\in[0,\bar p]$ and $p^1\leq p^2$ (respectively, $p^1\lneq p^2$).
\end{definition}

\vspace{.2cm}
\begin{lemma}\label{lem.optimum}
The polyhedron~\eqref{eq.d} has the following properties:
\begin{enumerate}
\item a \emph{maximal} point $p^*\in\calD$ exists that dominates all other points $p^*\geq p$, $\forall p\in\calD$;
\item $p^*$ is a global minimizer in the optimization problem
\beq\label{eq:clearing-opt-general}
\min_x\, F(x)\quad \text{subject to}\quad x\in\calD
\eeq
whenever function $F:[0,\bar p]\to\Real{}$ is non-increasing. If $F$ is decreasing, then $p^*$ is a \emph{unique} minimizer in~\eqref{eq:clearing-opt-general}.
\item $p^*$ is a clearing vector for the financial network;
\item each strongly connected component of graph $\calG[A]$, which is a sink (i.e., no arcs leave it), contains at least one node $i$ such that $p_i^*=\bar p_i$;
\item $p^*$ is the only clearing vector that enjoys property 4.
\end{enumerate}
\end{lemma}
\begin{remark}\rm
Statement~(2) is a refinement of~\cite[Lemma~4]{Eisenberg2001} stating that if $F$ is strictly decreasing, then every solution to~\eqref{eq:clearing-opt-general} is a clearing vector. Notice that formally the result in~\cite{Eisenberg2001} does not imply the uniqueness of a minimizer.
\end{remark}

\begin{figure}
  \centering
  \includegraphics[width=0.8\columnwidth]{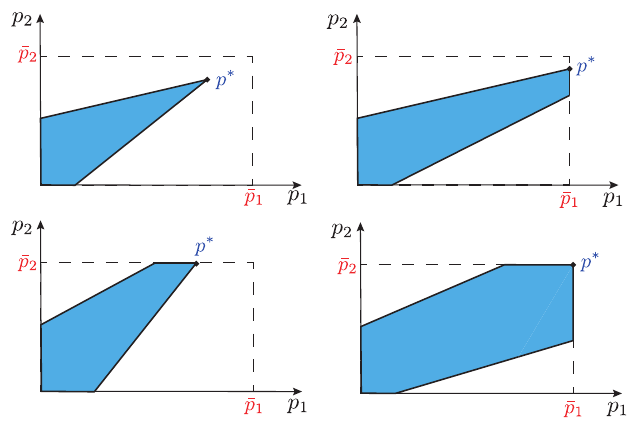}
  \caption{Examples of  possible shapes of  $\calD$ and the dominant clearing vector in the case of two nodes.}\label{fig.dominant}
\end{figure}

The clearing vector $p^*$ from Lemma~\ref{lem.optimum} is henceforth  referred to as the \emph{dominant} clearing vector, because it dominates all elements of $\calD$ (e.g., all possible clearing vectors). Figure~\ref{fig.dominant} illustrates possible structures of polyhedron $\calD$ and the location of $p^*$ in the case of $n=2$.

Introducing the vector of losses (the discrepancies between nominal and actual assets)
\[
\ell(p)=\bar\phi\ap{in}  - \phi\ap{in}=(\bar c-c)+A^{\top}(\bar p-p)\geq 0,
\]
statement~(2) of Lemma~\ref{lem.optimum} implies that a clearing vector can be found by solving a convex optimization problem such as, e.g., the convex QP
\bea
 \min_{p\in\Real{n}} & \|\ell(p)\|_2^2 \label{eq_clearingopt2}\\
 \mbox{s.t.:} & \bar p \geq p\geq 0 \nonumber \\
 & c + A\tran p  \geq p, \nonumber
 \eea
 where $\phi\ap{in} = c + A\tran p$,
 $ \bar \phi\ap{in} = \bar c + A\tran\bar p$,
or the conventional LP as follows
\bea
 \min_{p\in\Real{n}} & \one^{\top}\ell(p)  \label{eq_clearingopt_LP}\\
 \mbox{s.t.:} & \bar p \geq p\geq 0 \nonumber \\
 & c + A\tran p  \geq p. \nonumber
 \eea
Indeed, cost functions~\eqref{eq_clearingopt2} and~\eqref{eq_clearingopt_LP} are decreasing on $[0,\bar p]$. To prove this,
consider $p^1,p^2\in[0,\bar p]$ such that $p^1\lneq p^2$. Then, $\one^{\top}(p^2-p^1)=\sum_{i\in\calV}(p_i^2-p_i^1)>0$.
Recalling that $A$ is a stochastic matrix, it is now obvious that~\eqref{eq_clearingopt_LP} is strict decreasing:
\begin{equation}\label{eq.aux0}
\one^{\top}(\ell(p^1)-\ell(p^2))=\one^{\top}A^{\top}(p^2-p^1)=\one^{\top}(p^2-p^1)>0.
\end{equation}
To prove that~\eqref{eq_clearingopt2} is strict decreasing, notice that the vector
$\ell(p^1)-\ell(p^2)=A^{\top}(p^2-p^1)$ is nonnegative, furthermore, at least one its component is positive due to~\eqref{eq.aux0}.
On the other hand, $\ell(p^1),\ell(p^2)$ are also nonnegative, because $p^1,p^2\leq\bar p$. This implies that
\[
\ell(p^1)_i\geq\ell(p^2)_i\geq 0\quad\forall i\in\calV,
\]
and at least one of these inequalities is strict. Hence, $\|\ell(p^1)\|_2^2>\|\ell(p^2)\|_2^2$.

According to Lemma~\ref{lem.optimum}, therefore, a unique minimizer $p=p^*$ exists in the optimization problems~\eqref{eq_clearingopt2} and~\eqref{eq_clearingopt_LP}. Notice that usually uniqueness of a solution is guaranteed only for \emph{strictly} convex function, whereas the linear function~\eqref{eq_clearingopt_LP} is not strictly convex, and the function~\eqref{eq_clearingopt2} fails to be strictly convex if $\det A=0$.

%It should be noted that the optimization problems to find the dominant clearing vector admit various reformulations; some of them were used %in~\cite{LiuStaum:2010} to study the sensitivity of the clearing vector to system's parameters.

Although the dominant clearing vector $p^*$ depends on matrix $A$ and vectors $c,\bar p$, and the closed-form analytic expressions for its elements are not available, some properties of this clearing vector in fact do not depend on the liability matrix and are determined only by the topology of graph $\calG$. In particular, the following lemma can be proved
(see the appendix Section~\ref{sec:appendix} for its proof), which establishes the criterion for positivity of the dominant clearing vector's elements.

%\vspace{.2cm}
\begin{lemma}\label{lem.positivity}
The element $p_i^*$ of the dominant clearing vector is positive if and only if $i$ is not a sink node ($\bar p_i>0$) and at least one of the following conditions holds:
\begin{enumerate}
\item $i$ has outside assets, that is, $c_i>0$;
\item $i$ is reachable from some node $j\ne i$ with $c_j>0$;
\item the strongly connected component of graph $\calG$ to which $i$ belongs is a sink component.
\end{enumerate}
\end{lemma}
The result of Lemma~\ref{lem.positivity} proves to be useful in \emph{dynamic} (multi-step) models of interbank clearing that have been recently examined in~\cite{Calafiore2022,Calafiore2023}.

\subsection{Uniqueness of the clearing vector: a sufficient condition}\label{subsec.uniq-suffic}

In this subsection, we offer a \emph{sufficient} condition ensuring that the dominant clearing vector $p^*$ is the unique clearing vector. In fact, we show that some elements of the clearing vector are always determined uniquely.
We start by introducing some auxiliary notation. Let $C^+\doteq\{i:c_i>0\}$ stand for the set of nodes that receive nonzero outside assets,
 and let $S\doteq\{i:\bar p_i=0\}$ stand for the set of \emph{sink} nodes, who owe no liability payments. We introduce the set
\beq\label{eq.i0}
I_0\doteq C^+\cup S=\{i: c_i>0\lor \bar p_i=0\}.
\eeq
The following lemma, whose proof is given in the appendix Section~\ref{sec:appendix}, establishes a sufficient condition for uniqueness of the clearing vector and generalizes results from~\cite{Eisenberg2001,Glasserman2016,ElBitar2017}.

\begin{lemma}\label{lem.suffic}
Let $I_0'\supseteq I_0$ stand for the set of all nodes in the graph $\calG$, from where $I_0$ can be reached. Then, for \emph{every} clearing vector $p$ we have
$
p_i=p_i^*\quad\forall i\in I_0'.
$
In particular, if set $I_0$ is \emph{globally} reachable in the graph of a financial network $\calG$, then the dominant clearing vector $p^*$ is the only clearing vector corresponding to  the vector of outside assets $c$.
\end{lemma}
\vspace{.2cm}

\begin{remark}\rm
Since each path in the graph ends in one of the \emph{sink} components, it can be easily proved that {the uniqueness condition from the second part} of Lemma~\ref{lem.suffic} admits the following equivalent reformulation: each strongly connected component of graph $\calG$ being a sink is either trivial (contains only one node) or contains node $i$ such that $c_i>0$. In this form, the uniqueness criterion from Lemma~\ref{lem.suffic} becomes a special case of~\cite[Theorem~4.5]{Herings2021}. The latter theorem  deals with more general clearing systems, where the pro-rata rule is replaced by a more general division rule.
\end{remark}

\subsection{Uniqueness of the clearing vector: the general case}\label{subsec.uniq-gen}

In this subsection we derive two criteria of the clearing vector's uniqueness. The first of them (Theorem~\ref{thm.unique}) assumes that the dominant clearing vector $p^*$ has been found.
In this situation, we are able not only to check the uniqueness of the clearing vector, but also to describe the whole polytope of the clearing vectors without preprocessing the graph, e.g., computing the structure of its strongly connected components which is prerequisite for algorithms from~\cite{Hurd2016,Massai2022}.
On the other hand, if one is interested only in the uniqueness of a clearing vector, and the structure of the graph is known,
a simpler graph-theoretic criterion can be used (Theorem~\ref{thm.unique1}) that does not require knowledge of $p^*$.
{It should also be noticed that, unlike the initial results from~\cite{Eisenberg2001,Glasserman2016}, the uniqueness criteria from Theorem~\ref{thm.unique} and Theorem~\ref{thm.unique1} are not only sufficient, but also necessary.}

Assume that the condition in Lemma~\ref{lem.suffic} does not hold, that is, $I_0'\ne\calV$.  The banks corresponding to nodes from $\calV_1\doteq\calV\setminus I_0'$ do not have outside assets ($\calV_1\cap C^+=\emptyset$) and do not pay to nodes from $I_0'$ (otherwise, a chain of liability from them to $I_0$ would exist). Hence, matrix $A^1\doteq (a_{ij})_{i,j\in\calV_1}$ is stochastic.
At the same time, nodes from $I_0'$ \emph{can} have liability payments to nodes from $\calV_1$, which  depend only on the dominant vector $p^*$ and constitute the vector
\[
c^{(1)}\doteq (c^{(1)}_i)_{i\in\calV_1},\quad c^{(1)}_i=c^{(1)}_i(p^*)\doteq\sum_{k\in I_0'}a_{ki}p_k^*,\quad i\in\calV_1.
\]
We can now apply Lemma~\ref{lem.suffic} to a reduced financial network $\calG_1$ with node set $\calV_1$, normalized payment matrix $A_1$ and vector of external assets $c^{(1)}$. Introducing the set\footnote{Notice that unlike $I_0$,  set $I_1$ contains no sink nodes, because all sink nodes of the graph $\calG$ belong to $I_0$.}
\[
I_1=\{i\in\calV_1:c_i^{(1)}>0\}
\]
and denoting $I_1'\supseteq I_1$ all nodes from which set $I_1$ is reachable (banks that are connected by chains of liability to nodes from $I_1$), Lemma~\ref{lem.suffic} ensures that the elements of the reduced network's clearing vector $p_i,i\in I_1'$, are determined \emph{uniquely}. The definition~\eqref{eq:clearing-2} entails that if $p$ is a clearing vector for the original network, then its subvector $p^1=(p_i)_{i\in\calV_1}$ is a clearing vector for the reduced network $\calG_1$. This also applies to $p^*$. Lemma~\ref{lem.suffic} entails now that for each clearing vector $p$ (in the original network) one has $p_i=p_i^*\,\forall i\in I_1'$.

If $I_0'\cup I_1'=\calV$, we have uniqueness of the clearing vector. Otherwise, we have a group of banks $\calV_2=\calV\setminus(I_0'\cup I_1')$ that are not in debt to the nodes from $I_1'$ and $I_0'$, however, they can receive liability payments from the group $I_1'$. For group $\calV_2$, these payments may be treated as outside assets. Let
\[
c^{(2)}\doteq (c_i^{(2)})_{i\in\calV_2},\quad c_i^{(2)}\doteq\sum_{k\in I_1}a_{ki}p_i^*.
\]
If the set $I_2=\{i\in\calV_2:c_i^{(2)}>0\}$ is non-empty, one can consider the set $I_2'\supseteq I_2$
of all nodes from where $I_2$ can be reached. Lemma~1 implies that the elements $p_i,\,,i\in I_2'$ of the clearing vector are uniquely determined: $p_i=p_i^*\,\forall i\in I_2'$.

We arrive at the following iterative procedure, which allows to test the clearing vector's uniqueness (and, in fact, even to find the whole set of clearing vectors).
\begin{algorithm}[ht]\caption{Clearing vector's uniqueness test.}\label{alg.1}
\begin{algorithmic}
\State \textbf{Initialization.} Compute the dominant clearing vector $p^*$ (e.g., by solving the LP~\eqref{eq_clearingopt2}). Set $q\gets 0$,
$I_0\gets C^+\cup S=\{i:c_i>0\lor \bar p_i=0\}$. Find the set $I_0'\supseteq I_0$ of all nodes, from which $I_0$ is reachable.
\Repeat
    \State 1) $q\gets q+1$;
    \State 2) $\calV_q\gets\calV\setminus(I_0'\cup I_1\ldots\cup I_{q-1}')$;
    \State 3) compute the vector of payments from $I_{q-1}'$ to $\calV_q$
    \[
    c^{(q)}=(c^{(q)}_{i})_{i\in\calV_q},\quad c^{(q)}_{i}\doteq\sum_{k\in I_{q-1}'}a_{ki}p_k^*\quad\forall i\in\calV_q;
    \]
    \State 4) find the set $I_q=\{i\in\calV_q:c^{(q)}_i>0\}$;
    \State 5) find the set $I_{q}'\supseteq I_q$ of nodes from $\calV_q$, from where $I_q$ can be reached in $\calG$.
\Until{$\calV_q=\emptyset$ or $c^{(q)}=0$.}
\end{algorithmic}
\end{algorithm}

\begin{theorem}\label{thm.unique}
Algorithm~\ref{alg.1} stops after a finite number of steps $s\geq 0$. The elements of a clearing vector, corresponding to indices $i\in I_0'\cup I_1'\cup\ldots\cup I_s'$, are uniquely determined: $p_i=p_i^*$. The clearing vector is unique if and only if $\calV_s=\emptyset$, otherwise, there are infinitely many clearing vectors. Precisely, $p$ is a clearing vector if and only if
\beq\label{eq.all-clearing-vec}
p_i=\begin{cases}
p_i^*,\quad i\in I_0'\cup I_1\ldots\cup I_s',\\
\xi_i,\quad i\in\calV_s,
\end{cases}
\eeq
where $\xi\in\Real{\calV_s}$ is an arbitrary vector satisfying the constraints
\beq\label{eq.subvector-constraint}
B^{\top}\xi=\xi,\quad 0\leq\xi_i\leq\bar p_i,\;\forall i\in\calV_s,\quad B\doteq (a_{ij})_{i,j\in\calV_s}.
\eeq
\end{theorem}
A proof of Theorem~\ref{thm.unique} is offered in the appendix Section~\ref{sec:appendix}.
Theorem~\ref{thm.unique} entails the following technical proposition, which is of independent interest.
 The \emph{net} vector (or the vector of \emph{equities}) corresponding to the clearing vector $p$ is
\beq\label{eq.equity}
\zeta\doteq \phi\ap{in}-\phi\ap{out}=c+A\tran p-p=\max(c+A\tran p-\bar p,0)\geq 0.
\eeq
The component $\zeta_i$ represents the \emph{net worth} of node $i$, that is, the difference between the in- and the out-flows at that node.
\begin{corollary}\cite{Eisenberg2001,Kabanov2018}\label{cor.equity}
The vector of equities $\zeta$ is the same for all possible clearing vectors $p$.
\end{corollary}

The proof of Corollary~\ref{cor.equity} is straightforward from the representation~\eqref{eq.all-clearing-vec}. Indeed, the equities of banks from $\calV_s$ are all zeros, whereas the remaining equities depend only on the subvector $(p_i)_{i\in I_0'\cup\ldots\cup I_s'}$, which is uniquely determined by $p^*=p^*(c,A,\bar p)$.\qed

Notice that although the subvector $\xi$ in~\eqref{eq.all-clearing-vec} is defined non-uniquely, some of its elements are in fact uniquely determined due to Lemma~\ref{lem.positivity}. As we know, $\xi_i=p_i^*=0$ whenever $i$ does not belong to a sink component and is not reachable from $C^+$. Combining Theorem~\ref{thm.unique} with  Lemma~\ref{lem.positivity}, we can establish an alternative uniqueness criterion, which \emph{does not} require knowledge of the vector $p^*$.

\vspace{.2cm}
\begin{theorem}\label{thm.unique1}
 The following  two conditions are equivalent:
  \begin{enumerate}[label=\roman*)]  \item the clearing vector is unique; %(and equals $p^*$);
  \item each non-trivial sink component of $\calG$ either contains a node from $C^+$ or is reachable from $C^+$.
  \end{enumerate}
\end{theorem}

\noindent
A proof of Theorem~\ref{thm.unique1} is given in the appendix Section~\ref{sec:appendix}.

\subsection{Discussion: Lemma~\ref{lem.suffic} and Theorems~\ref{thm.unique} and~\ref{thm.unique1} vs. previously known results}

Comparing Theorem~\ref{thm.unique} to previously known results, several relevant features can be noticed. To find the dominant clearing vector $p^*$, one needs to know the matrix $\bar P$ (or, equivalently, $\bar p$ and $A$) and the vector of external assets $c\geq 0$. As we have already discussed, efficient algorithms do exist for computing this vector (e.g., by solving an LP or a convex QP). The computation of the sets $I_s,I_s',\calV_s$ requires in fact very limited information. We do not need to compute the vectors $c^{(q)}$, it suffices to know which of their components are positive. A closer look at our algorithm shows that this depends on the set $C^+$ (nodes who receive external assets)
and the topology of graph $\calG$ (which determines, in view of Lemma~\ref{lem.positivity}, which elements of $p^*$ are positive). Also, we give an explicit and simple description of \emph{all} clearing vectors.

\subsubsection{Criteria from~\cite{Eisenberg2001} and~\cite{Glasserman2016}}\label{subsubsect.discuss.1}

{The second part of Lemma~\ref{lem.suffic} (the uniqueness of a clearing vector if $I_0$ is globally reachable) covers two well-known criteria of the clearing vector's uniqueness.}

The uniqueness criterion from~\cite[Theorem~2]{Eisenberg2001} states that the clearing vector is unique if the set $C^+$ is globally reachable. {In~\cite{Eisenberg2001}, this ``regularity'' property is formulated as follows: the ``risk orbit'' of each node $i$ (that is, the set of nodes reachable from $i$) is a surplus set, that is, the orbit contains a node $j$ with $c_j>0$. Notice that formally this definition of regularity does not apply to a network with a sink node, whose ``risk orbit'' is empty. It can be proved that Theorem~2 in~\cite{Eisenberg2001} retains its validity if one formally considers the degenerate orbit of a sink node as a surplus set.} An advantage of the approach
developed in~\cite{Eisenberg2001} is the possibility to generalize the proof (with some variations) to some advanced models of systemic risk~\cite{Kabanov2018,Elsinger2006}.

{The criterion from~\cite{Glasserman2016} guarantees that the clearing vector is unique if all nodes are connected to the external sector by chains of liabilities. Introducing the fictitious sink node, this can be formulated as follows: the set $S$ of sink nodes contains only one node, which is globally reachable.}

{Note that the first part of Lemma~\ref{lem.suffic}, although it follows from the second part, is not easily available in the literature. This statement plays the central role in our uniqueness criteria (Theorems~\ref{thm.unique} and~\ref{thm.unique1}).} Unlike our Theorem~\ref{thm.unique}, the aforementioned results from~\cite{Eisenberg2001,Kabanov2018,Glasserman2016} do not allow to parameterize all clearing vectors when the uniqueness fails to hold.

\subsubsection{An extension of the Eisenberg-Noe model: the negative inflow vector}
One of the limitations of the Eisenberg-Noe model is the assumption that banks cannot become \emph{insolvent}, that is, the equity vector~\eqref{eq:in-out-flows} is always nonnegative. As discussed in~\cite{Elsinger2006,Hurd2016} the latter assumption can be violated in practice, because some payments to the external sector (e.g., the bank's operation costs)
are more senior than debts to other banks and cannot be reduced, even if the external assets drop. Hence, a generalization of the Eisenberg-Noe model has been introduced in~\cite{Elsinger2006} where the vector $c$ can have both positive and negative components $c_i\in\mathbb{R}$, which stand for the difference between outside asset of bank $i$ and its external liability.
In this situation, the definition of the clearing vector has to be modified to
\beq\label{eq.elsinger}
p=\min\left(\bar p,\max(A\tran p+c,0)\right),
\eeq
which is equivalent to~\eqref{eq:clearing-2} if $c\geq 0$.

The existence of clearing vectors in the case of sign-indefinite $c$ can be derived from the Knaster-Tarski theorem, and the standard fixed-point iteration delivers one of the clearing vectors~\cite{Elsinger2006}. However, the uniqueness problem has been addressed quite recently. For strongly connected networks, the uniqueness was established in~\cite{Acemoglu2015} and~\cite{Ren2016}. For the general situation, two different uniqueness criteria have been proposed in~\cite[Theorem~2.2]{Hurd2016} and~\cite[Theorem~2]{Massai2022}, which also provide a description of all possible clearing vectors.
In the situation when $c\geq 0$, however, the cited criteria appear to be superfluous and  more computationally demanding than our results.
Notice that our Theorem~\ref{thm.unique} requires finding the maximal
clearing vector $p^*$ (which, in the case $c\geq 0$, is obtained by solving a convex QP or an LP), however, this theorem does
not require knowledge of the structure of the strongly connected components of the graph (or, equivalently, of an irreducible decomposition of
matrix $A$). On the other hand, if the structure of strongly connected sink components is known, Theorem~\ref{thm.unique1} allows
to test the clearing vector's uniqueness without finding the clearing vector itself. The general
results from~\cite{Hurd2016,Massai2022} exploit the irreducible decomposition of matrix $A$ (whose computation is a self-standing
non-trivial procedure for a large-scale network), and do not allow to test the uniqueness without computing the complete clearing vector.
Notice that the procedure of computing the clearing vector (even in the case of its uniqueness) is not made explicit in~\cite{Hurd2016} (one of the option is the modified fictitious default algorithm~\cite{Elsinger2006}).
The procedure of clearing vector computation in~\cite{Massai2022} is explicit, however, it requires finding the left and right Perron-Frobenius eigenvectors for each irreducible block, which is not needed in our results.

\subsection{Examples}
\begin{example}\rm
We first illustrate our uniqueness criterion on the synthetic financial network shown in Figure~\ref{fig.3comp}. This network includes one sink node, representing the external sector, and three non-trivial strongly connected components (i.e., cliques).

The strongly connected component 1 stands for a group of banks that works directly with the external sector, using also some financial services of banks from group 3 and thus owing some liabilities to them.
The banks from group 2, constituting another strongly connected component, owe debts to banks from groups 1 and 3; these banks thus also have indirect liability to the external sector (through banks from group 1).
The group of banks 3 constitutes a non-trivial \emph{sink} component, owing neither direct nor even indirect liability to the financial sector. Notice that the presence of this component makes the criterion from~\cite{Glasserman2016} inapplicable: the set $S$ of sink nodes is not globally reachable.

Denoting the vector of outside assets corresponding to group $i$ by $c_{[i]}\geq 0$, the following situation are possible:
\begin{itemize}
\item If $c_{[3]}\ne 0$, we are in the situation of~\cite[Theorem~2]{Eisenberg2001}: the set $C^+$ is reachable from any node but for the sink node (as has been remarked, the result from~\cite{Eisenberg2001} retains its validity in presence of sinks).
\item If $c_{[3]}=0$, then the result from~\cite{Eisenberg2001} becomes inapplicable (group 3 is not connected to $C^+$). Nevertheless, Theorem~\ref{thm.unique1} ensures the uniqueness of the clearing vector if $c_{[1]}\ne 0$ or $c_{[2]}\ne 0$, because component $3$ is reachable from set $C^+$.
\item Finally, if $c=0$, the clearing vector is not unique. However, the subvector corresponding to the union of groups 1 and 2 is unique: $p_{[1]}=p_{[1]}^*$ and $p_{[2]}=p_{[2]}^*$ for each clearing vector $p$ (here the subscript $[i]$ has the same meaning as for subvectors $c_{[i]}$). This is entailed by the first part of Lemma~\ref{lem.suffic}, because groups $1$ and $2$ are connected to the unique sink node (set $S$).
\end{itemize}
\begin{figure*}[t]
\centering
\includegraphics[width=0.75\columnwidth]{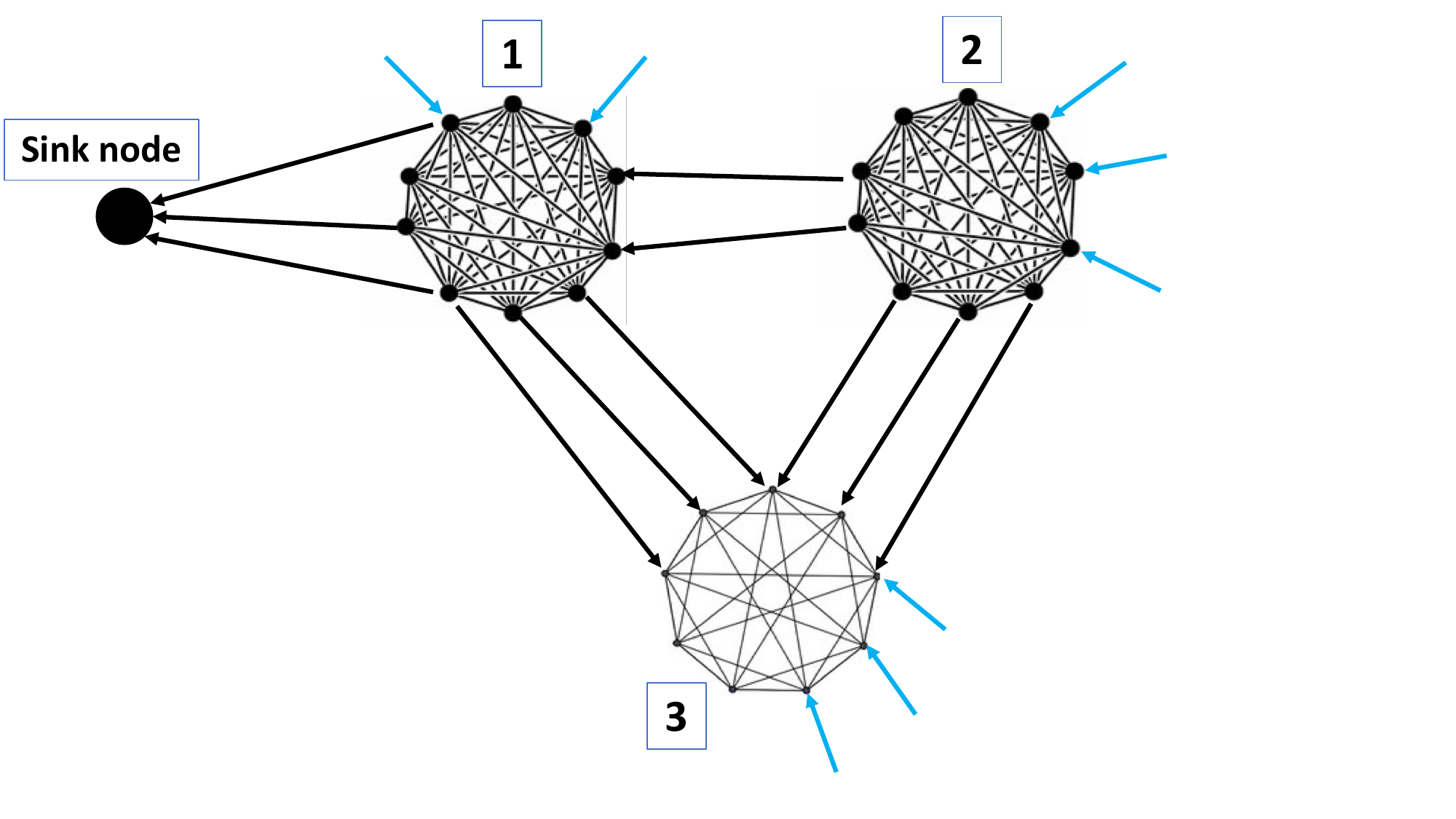}
\caption{Network with a sink node and three strongly connected components.}\label{fig.3comp}
\end{figure*}
\end{example}
\begin{example}\rm
Next, we illustrate Algorithm~\ref{alg.1} by considering the network in Figure~\ref{fig.netw-big}, which displays a network with $n=15$ nodes
that contains only one sink node ($S=\{0\}$) and three nodes with outside assets ($C^+=\{1,2,3$\}), which together form
 the set $I_0$.
Lemma~\ref{lem.positivity} entails that in this situation $p_i^*>0$, $\forall i\ne 0,9$ (notice that nodes 13-15 constitute a sink component, satisfying thus condition 3) from Lemma~\ref{lem.positivity}).
The set $I_0'$  contains $I_0$ and two nodes $4,5$ that have liabilities towards nodes $0$ and $2$.
The set $I_1$ contains nodes that have no liability to $I_0'$, however, should receive  payments from 4 and 5.
Hence,  $c^{(1)}_6,c^{(1)}_7,c^{(1)}_8>0$. The nodes $6,7,8$ constitute the set $I_1$; the set $I_1'$ is obtained by adding node $9$ who has  liability to one of them.
On the next iteration of the algorithm, one computes the sets $I_2=\{10,11\}$ and $I_2'=\{12\}\cup I_2$.
The next vector $c^{(3)}$ will be zero, because the remaining nodes of the graph constitute an isolated group.
Hence, the clearing vector is not unique, but for each clearing vector $p$ one has $p_i=p_i^*$, $\forall i=0,\ldots,12$.
\begin{figure*}[t]
\centering
\includegraphics[width=0.75\columnwidth]{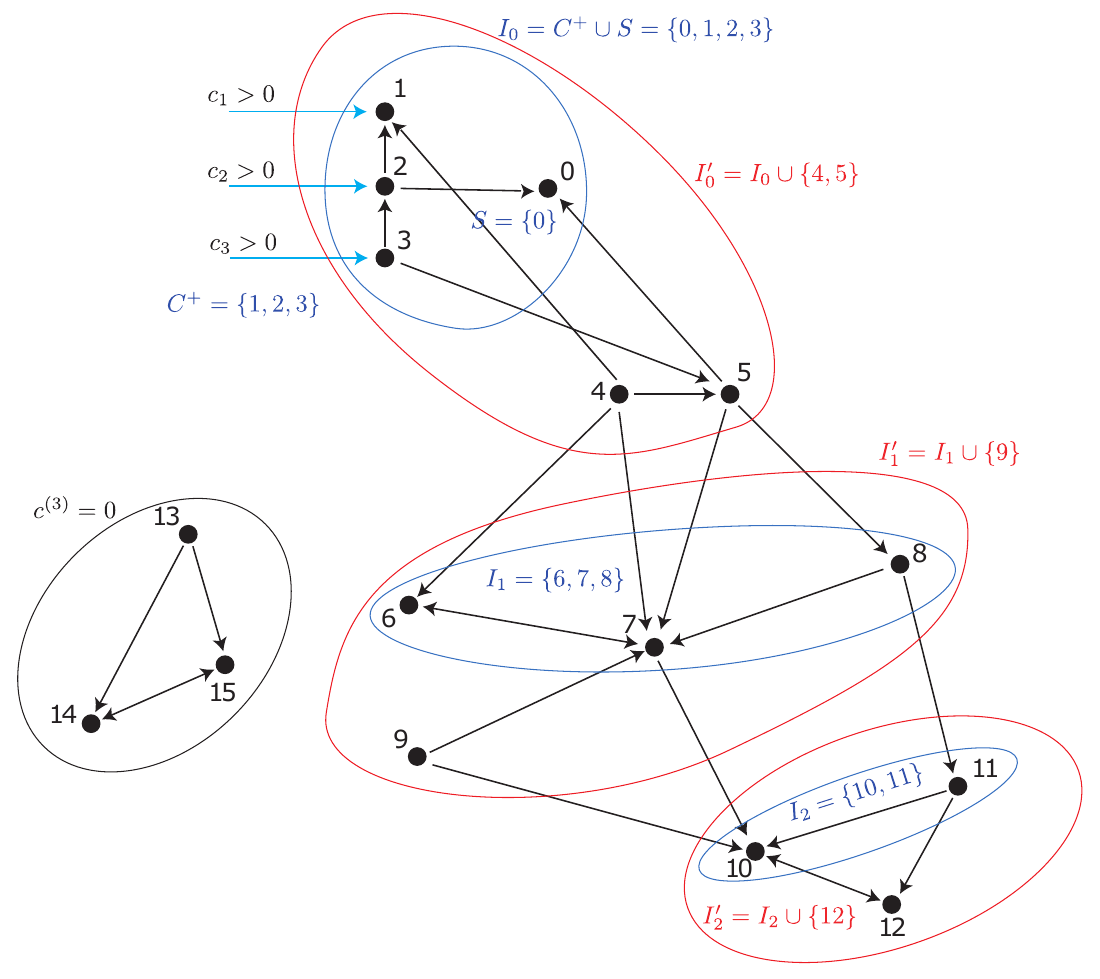}
\caption{The sets $I_s$ (encircled by blue lines) and $I_s'$ (encircled by red lines) for a special financial network with $n=15$ nodes with $\bar p_i>0$, $\forall i\ne 0$.}\label{fig.netw-big}
\end{figure*}

Removing the isolated strongly connected component $\{13,14,15\}$ from the graph in Figure~\ref{fig.netw-big}, one obtains a financial network that has only one clearing vector $(p_i^*)_{i=1}^{12}$ yet does not satisfy the sufficient condition from Lemma~\ref{lem.suffic}:
set $C^+\cup S$ cannot be reached from any of the nodes $6,\ldots,12$.
None of the criteria from~\cite{Eisenberg2001} and~\cite{Glasserman2016} guarantees uniqueness in this situation.
\end{example}

\begin{remark} \rm
In the situations considered in Examples 1 and 2 the outside assets of several banks drop simultaneously ($c_i=0$ for several nodes $i$).
This situation is common in actual financial networks where the banks are connected not only by mutual liabilities, but also by mutual exposures to external risky assets~\cite{Allen2012,Elliott2014}. A bankruptcy of some business or a forced sale of some asset at a low price in the case of the market's turbulence~\cite{Cifuentes2005} has a negative effect on all shareholders, as illustrated by the housing market crash in 2008. Such ``systematic'' shocks are introduced in problems of portfolio compression~\cite{Amini2023} and have to be considered in any realistic stress-test simulation.
\end{remark}

\section{Lifting the pro-rata condition}\label{sec:no_prorata}
In this section we propose a new clearing mechanism that does not hinge upon the proportional payments rule.
While such new clearing mechanism is not yet in operation, we here advocate the idea that a more transparent knowledge of the inter-bank liability matrix, and a centralized authority in charge of defining the clearing payments in case of default, may lead to important systemic and social benefits due to the avoidance or reduction of catastrophic default events.

Without the pro-rata constraint in place, we have potentially more freedom to appropriately select a payment matrix $P$ that satisfies the conditions of
Definition~\ref{def.matrix}. These conditions do not identify the payment matrix univocally. Indeed,  as we discussed extensively in the previous sections, one matrix that satisfies the definition is the one resulting from the pro-rata rule, but this may not be the only one, in general.
This fact gives us some degrees of freedom that shall be used towards seeking clearing matrices who promote a virtuous behavior of the whole network. As it will be shown in Section~\ref{sec:experiments}, the proposed non-proportional clearing mechanism
not only visibly reduces the overall  total imbalance between the nominal and actual payments, but it also helps isolating defaults and preventing their cascaded spread over the network.

%There are in fact several reasons to consider clearing matrices different from the pro-rata matrix. The first reason is that, as discussed %in~\cite{csoka2018decentralized,Glasserman2016}, the creditors of different standings may have different priorities, and \emph{the principles of %proportionality and priority play an important role in bankruptcy law across the globe}~\cite{csoka2018decentralized}. Another reason is %that the relaxation of the pro-rata constraint can visibly reduce a systemic loss function~\cite{Glasserman2016}, that is, the total imbalance %between the nominal and actual payments (see Section~\ref{sec:experiments}).

\subsection{Optimal clearings}
We next present an optimization-based approach for determining a clearing matrix that minimizes an appropriate system-level loss function.
Recalling that $[0,\bar P]=\{P\in\Real{n\times n}:0\leq P\leq\bar P\}$, we consider the convex polyhedron in the space of matrices %$
\[
\calD_{n\times n}=[0,\bar P]\cap\{P\in\Real{n\times n}:c+P\tran\one\geq P\one\}.
\]
We call a function $F:[0,\bar P]\to\Real{}$ \emph{decreasing} if $F(P_1)>F(P_2)$ whenever $P_1\lneq  P_2$. For any such function, consider the optimization problem
\beq\label{eq:clearing-opt-general-matr}
\min F(P)\quad \text{subject to}\quad P\in\calD_{n\times n}.
\eeq
The following result holds; see Section~\ref{sec:appendix} in the Appendix for a proof.
\begin{lemma}\label{lem:opt_no_prorata}
For any decreasing function $F:[0,\bar P]\to\Real{}$, every local minimizer in
problem~\eqref{eq:clearing-opt-general-matr} is a clearing matrix (as defined in Definition~\ref{def.matrix}).
\end{lemma}

Notice that Lemma~\ref{lem:opt_no_prorata} does not require the function to be continuous. For
a continuous function, the global minimum always exists due to compactness of $\mathcal{D}_{n\times n}$.
Two examples of functions that are continuous and decreasing on $[0,\bar P]$ are
\begin{eqnarray}
F_1(P) \doteq \|\bar P  - P\|_1 &= &\sum_{i,j=1}^n |\bar P_{ij}- P_{ij}|  ,\label{eq.loss-lin} \\
F_2(P) \doteq
\|\bar P  - P\|_2^2 &= &\sum_{i,j=1}^n (\bar P_{ij}- P_{ij})^2 .\label{eq.loss-quad}
\end{eqnarray}
Both these functions provide an aggregated, system-level, measure of the impact of defaults on the deviation of payments from their nominal liabilities.
In particular, minimization of the performance index~\eqref{eq.loss-lin} is equivalent to the  problem
\bea
F_1^* =  \min_{P\in\Real{n,n}} & \|\bar P  - P\|_1  \label{eq_clearingopt_free_inflow1}\\
 \mbox{s.t.:} & \bar P \geq P \geq 0 \nonumber \\
 & c + P\tran \one - P \one  \geq 0,\nonumber
 \eea
 which can be readily recast in the form of an LP, whereas minimization of~\eqref{eq.loss-quad}  leads to the  convex QP problem
\bea
F_2^* =  \min_{P\in\Real{n,n}} & \|\bar P  - P\|_2^2  \label{eq_clearingopt_free_inflow}\\
 \mbox{s.t.:} & \bar P \geq P \geq 0 \nonumber \\
 & c + P\tran \one - P \one  \geq 0.\nonumber
 \eea
 While both \eqref{eq_clearingopt_free_inflow1} and \eqref{eq_clearingopt_free_inflow}  are valid ways for obtaining a clearing matrix $P$, from the financial point of view they have different characteristics.
Notice in fact that the objective in \eqref{eq_clearingopt_free_inflow} is strongly convex, so the optimal clearing matrix resulting from
\eqref{eq_clearingopt_free_inflow} is unique. However, the $F_2$ criterion is expressed in ``squared currency,'' which may not have an immediate financial meaning. Further, this sum-of-squared losses function tends to avoid large individual losses,  at the expense of possibly having many small losses. This is a well-known behavior of squared losses, which may be critical  in the present context  since default is an on/off process: a node goes into default as soon as its balance goes negative, irrespective of how large is the negative balance.
Especially for this reason our focus in this work is on the $F_1$ criterion in  \eqref{eq_clearingopt_free_inflow1}. As it is widely known, this
$\ell_1$-norm criterion tends to be less sensitive to large residuals and promotes solutions that are {\em sparse}. Sparsity, in our context, is a very desirable property, since a sparse residual matrix $\bar P  - P$ means that many of its entries are zero, which in turn means that  many nodes  meet their liabilities  and thus do not default. We next focus on the $F_1$ criterion, which represents the total sum that is lost due to defaults (notice that in the case of no defaults we have $F_1^* = 0$).

\subsection{Non-uniqueness of clearing matrices}

In general, the set of clearing matrices defined via Definition~\ref{def.matrix} has more than one element.
Further, this set  has a non-trivial structure and fails to be a complete lattice\footnote{It should be noticed that some alternative definitions of clearing matrices are possible, under which the set of clearing matrices becomes a complete lattice, see, e.g., a very general construction from~\cite{Csoka2018} that deals with integer payments and allowing the banks to apply different bankruptcy policies.}.
Lemma~\ref{lem.optimum} does not retain its validity in the class of all admissible clearing matrices, in particular, there exists no maximal clearing matrix, as exemplified by the following example.

\begin{example}\rm
Consider a degenerate 3-node network whose nodes 1 and 2 are sinks ($\bar p_{1i}=0,\bar p_{2i}=0$, $ \forall i$) and receive no outside assets ($c_1=c_2=0$), whereas node 3 receives outside asset
$c_3>0$ and owes to both 1 and 2 ($\bar p_{31},\bar p_{32}>0$). A clearing matrix $P$ must then have the following structure
\beq\label{eq.aux2a}
P=\begin{bmatrix}
0 & 0 & 0\\
0 & 0 & 0\\
p_{31} & p_{32} & 0\end{bmatrix},\quad p_{31}\in [0,\bar p_{31}],\,p_{32}\in [0,\bar p_{32}].
\eeq
and equation~\eqref{eq:clearing-1} is equivalent to
\beq\label{eq.aux2}
p_{31}+p_{32}=\theta\doteq\min(\bar p_{31}+\bar p_{32},c_{3}).
\eeq
Considering two clearing matrices
\[
P_1=\begin{bmatrix}
0 & 0 & 0\\
0 & 0 & 0\\
\theta & 0 & 0\end{bmatrix}, \quad P_2=\begin{bmatrix}
0 & 0 & 0\\
0 & 0 & 0\\
0 & \theta & 0\end{bmatrix},
\]
one may easily notice that no matrix $P\geq \max(P_1,P_2)$ can satisfy~\eqref{eq.aux2}. In particular, the set of clearing matrices understood in the sense of Definition~\ref{def.matrix} has \emph{no} maximal element. It can be easily verified that every matrix~\eqref{eq.aux2a} that satisfies~\eqref{eq.aux2} is a \emph{global} minimizer of problem~\eqref{eq_clearingopt_free_inflow1}. Hence, in spite of the strict monotonicity of the cost function in~\eqref{eq_clearingopt_free_inflow1}, the optimal solution to~\eqref{eq_clearingopt_free_inflow1} is not unique.
Contrary, if we introduce the proportionality rule, the unique (thanks to Lemma~\ref{lem.suffic}) clearing vector and the corresponding clearing matrix, result to be, respectively
\[
p^*=(0,0,\theta)^{\top},\quad P^*=\begin{bmatrix}
0 & 0 & 0\\
0 & 0 & 0\\
\frac{\bar p_{31}\theta}{\bar p_{31}+\bar p_{32}} & \frac{\bar p_{32}\theta}{\bar p_{31}+\bar p_{32}} & 0\end{bmatrix}.
\]
\end{example}

\subsection{A two-stage approach for uniqueness}
\label{sec:two-stage}
Given a liability matrix $\bar P$ and an in-flow vector $c\geq 0$ our purpose is to establish a policy that permits  determining uniquely
a clearing matrix  that globally minimizes the system-level loss $F_1$.
Uniqueness is important since the clearing payments must be non controversial, and each player must be in condition to compute and agree univocally on the same payment values. As we discussed in the previous section, however, the LP problem in \eqref{eq_clearingopt_free_inflow1} may have multiple global optimal solutions,  in general.
The following lemma characterizes {\em all} optimal solutions to \eqref{eq_clearingopt_free_inflow1}.

\begin{lemma}\label{eq:lem_allsols}
Let $P^*$ be an optimal solution of problem~\eqref{eq_clearingopt_free_inflow1}. Then, the set of all optimal solutions to \eqref{eq_clearingopt_free_inflow1} is characterized as the polytope
\[
\calS^* \doteq \{\tilde P = P^*+\Delta :\; 0\leq  \tilde P \leq \bar P, \, c+\tilde P\tran \one -  \tilde P \one \geq 0,\, \one\tran \Delta \one = 0 \}.
\]
\end{lemma}
A proof of Lemma~\ref{eq:lem_allsols} is given in the appendix Section~\ref{sec:appendix}.

The next proposition illustrates how we can
break the tie and specify a two-stage rule that results in a unique optimal clearing policy.
\begin{proposition} \label{eq:prop_unique}
Let $P^*$ be an optimal solution of problem~\eqref{eq_clearingopt_free_inflow1}, and let $F_1^*$ be the corresponding globally optimal loss level. Let further $\Delta^*$ be such that
\bea
\Delta^* = \arg\min_{\Delta\in\Real{n,n}} & \|P^* + \Delta\|_2^2  \label{prob:delta_unique}\\
\mbox{s.t.:} & \bar P \geq  P^*+\Delta \geq 0 \nonumber \\
& c+(P^*+\Delta)\tran \one -  (P^*+\Delta) \one \geq 0 \nonumber \\
& \one\tran \Delta \one = 0. \nonumber
\eea
Then:
\bit
\item[(a)] $\tilde P^* \doteq P^* + \Delta^*$ is a clearing matrix;
\item[(b)] $\tilde P^*$ achieves the globally optimal loss level $F_1^*$, that is $F_1(\tilde P^*) = F_1^*$;
\item[(c)] $\tilde P^*$ is the unique  smallest-Euclidean-norm solution among all optimal solutions to \eqref{eq_clearingopt_free_inflow1}.
\eit
\end{proposition}
A proof of Proposition~\ref{eq:prop_unique} is given in the appendix Section~\ref{sec:appendix}.

We call this  a two-stage approach for determining the unique optimal clearing  matrix, since the method consists in solving two convex problems in sequence: first we solve the LP  \eqref{eq_clearingopt_free_inflow1} and find any optimal solution $P^*$. Next, we solve the convex QP \eqref{prob:delta_unique} and find its unique optimal $\Delta^*$. Finally, our optimal clearing matrix  of interest is $\tilde P^* = P^* + \Delta^*$,
which is the unique minimum norm solution among all possible optimal solutions to problem  \eqref{eq_clearingopt_free_inflow1}.

\section{Numerical Experiments}\label{sec:experiments}
We next evaluate numerically the systemic improvement of the optimal matrix clearings (i.e., those computed according to  Proposition~\ref{eq:prop_unique}) with respect to the standard pro-rated clearings (computed according to \eqref{eq_clearingopt_LP}).
The improvement is evaluated both in terms of reduction of the systemic loss $F_1$ and in terms of containment of default contagion, as expressed by the percentage of  defaulted nodes in the network.
First, we propose in Section~\ref{sec:schematicex} a simple schematic example, and then we perform
extensive randomized tests in Section~\ref{sec:ex:random} using synthetic random networks similar to ones proposed in~\cite{Nier2007}.

\subsection{An illustrative schematic example}\label{sec:schematicex}
 We consider a variation on the simplified network discussed in \cite{Glasserman2016}.
This network is composed of $n=5$ nodes (including one sink node), with liability matrix
\[
\bar P =  \left[\begin{array}{ccccc} 0 & 180 & 0 & 0 & 180\\ 0 & 0 & 100 & 0 & 100\\ 90 & 0 & 0 & 100 & 50\\ 150 & 0 & 0 & 0 & 150\\ 0 & 0 & 0 & 0 & 0 \end{array}\right].
\]
%where the last row refers to the  fictitious node.
%
Suppose there is a nominal  scenario where external cash flows are given as
\[
c = c\ped{nom} \doteq [121,\; 21,\; 150,\; 204,\; 0]\tran .
\]
It can be readily verified that in the nominal scenario all the nodes in the network remain solvent (no defaults), and the clearing payments coincide with the nominal liabilities.
Consider next a situation in which a  ``shock'' happens on the inflow at node $3$, so that this inflow reduces from $150$ to $130$, that is
  \[
c = c\ped{shock}\doteq  [121,\; 21,\; 130,\; 204,\; 0]\tran  .
\]

\begin{figure}[htb]
\centering
\includegraphics[width=1\columnwidth]{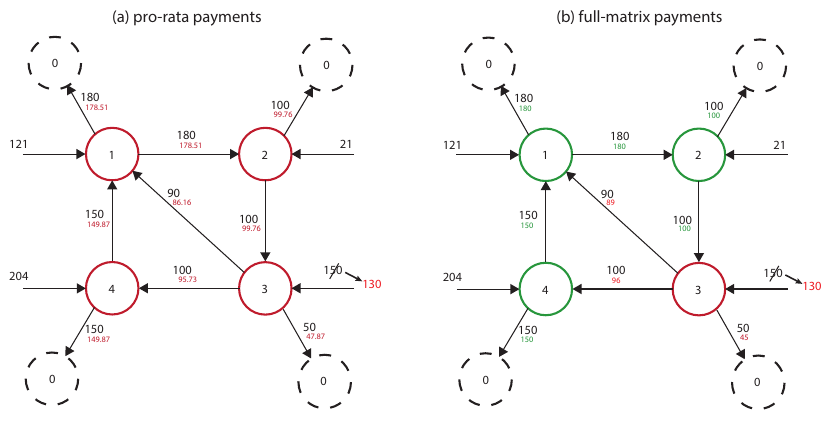}
\caption{A four-nodes network with a shock on the cash inflow at node 3. Panel (a) shows the pro-rated clearing payments, panel (b) shows the optimal matrix clearings.}\label{fig:glass_example}
\end{figure}

\noindent
Under the pro-rata rule, the (unique) clearing payments, resulting from the solution of  \eqref{eq_clearingopt_LP}, are shown in smaller font below the nominal liabilities in the left panel of Figure~\ref{fig:glass_example}: all nodes in the network default in a cascade fashion due to initial default of node $3$. The total defaulted amount, defined as the sum of all the unpaid liabilities is in this case $F_1=13.98$.

Then, we lifted the pro-rata rule, and we computed the unique clearing payments
according to Proposition~\ref{eq:prop_unique}. The results in this case are shown in the right panel of Figure~\ref{fig:glass_example}:
only node $3$ defaults, while all other nodes manage to pay their liabilities in full. Not only we reduced the  sum of all unpaid liabilities
 to $F_1=10$, but we also obtained {\em isolation} of the contagion, since the default is now
 limited to node $3$ and did not spread to other parts of the network.

\subsection{Random networks test} \label{sec:ex:random}
The random graphs used for simulations are constructed using a technique inspired by~\cite{Nier2007}.
The topology of the graph is given by the standard Erd\"os-Renyi ${\cal G}(n,p)$ graph.
The interbank liabilities $\bar{P}_{ij}$ for every edge $(i,j)$ of the random graph are then found by sampling from a uniform distribution $\bar{P}_{ij}\sim\mathcal{U}(0,P_{max})$, where $P_{max}$ is the maximum possible value of a single interbank payment. In the experiments we set $P_{max}=100$. Unlike~\cite{Nier2007}, the values $\bar P_{ij}$ can thus be heterogeneous.
%Also, we do not consider payments to external sector (deposits etc.): as it has been discussed, we can always get rid of them by introducing a fictitious node.

 Following~\cite{Nier2007}, we define the total amount  of the external assets  $E=\frac{\beta}{1-\beta}I$, where $I=\sum_{i,j=1}^n \bar P_{ij}$ is the total amount of the interbank liabilities
 and $\beta=E/(E+I)$ is a parameter representing the percentage of external assets in total assets at the system level; in our experiments $\beta=0.05$.
 The nominal asset vector $\bar c$ is then computed in two steps. First,  each bank is given the minimal value of external assets under which its balance sheet is equal to zero. At the second step, the remainder of the aggregated external assets is evenly distributed among all banks.

The financial shock is modeled by randomly choosing a subset of $n_s$ banks of the system and nullifying their external financial assets.

\subsubsection{A numerical study of the ``price'' of proportionality}

To evaluate the ``price'' of imposing the pro-rata rule, we consider the $F_1$ performance index presented in~\eqref{eq.loss-lin}. This is used as \emph{system loss} also in~\cite{Glasserman2016}.
%Obviously, this function is strictly monotone.
Its minimal value over all matrices obeying the pro-rata constraint~\eqref{eq:prorata} is
$F_1^{(\mathrm{pr})}=\sum_{i}(\bar p_i-p_i^*)$, where $p^*$ is the dominating clearing vector from Lemma~\ref{lem.optimum}, which is found by solving problem~\eqref{eq_clearingopt2}. Relaxing the pro-rata constraint, we use the two-stage approach presented in Section~\ref{sec:two-stage} to find the globally optimal clearing matrix $P^*$, resulting in the system loss $F_1^{(\mathrm{nopr})}=\|\bar P  - P^*\|_1$.
The price, or global effect,  of the pro-rata rule can thus be estimated by the following ratio
\[
G=\frac{F_1^{(\mathrm{pr})}-F_1^{(\mathrm{nopr})}}{F_1^{(\mathrm{nopr})}}\in [0,1].
\]
If $G=0$ (as, e.g., in Example~2, where all clearing matrices are optimal), the imposition of pro-rata constraint is ``gratuitous'' in the sense that it does not increase the aggregate system loss: $F_1^{(\mathrm{pr})}=F_1^{(\mathrm{nopr})}$.
The larger value $G$ we obtain, the more ``costful'' is the pro-rata restriction.
%\textcolor{red}{[GC: Above you talk about the cost of using pro-rata. The cost of using pro-rata should maybe be more properly defined as
%\[
%G=\frac{ F_1^{(\mathrm{pr})} - F_1^{(\mathrm{nopr})}}{F_1^{(\mathrm{nopr})} }\in [0, 1],
%\]
%that is, using the no-prorata cost as a baseline.
%This way, since the denominator is smaller than with the previous definition, this would give more ``dramatic'' numerical results. Can we redo Figure 4 using this slightly modified definition?
%Also, the caption  in  Figure 4 should be ``Cost of introducing the pro-rata rule'', and not ``gain'' obained by removing it...
%]}

\begin{figure}
  \centering
  \includegraphics[width=0.75\columnwidth]{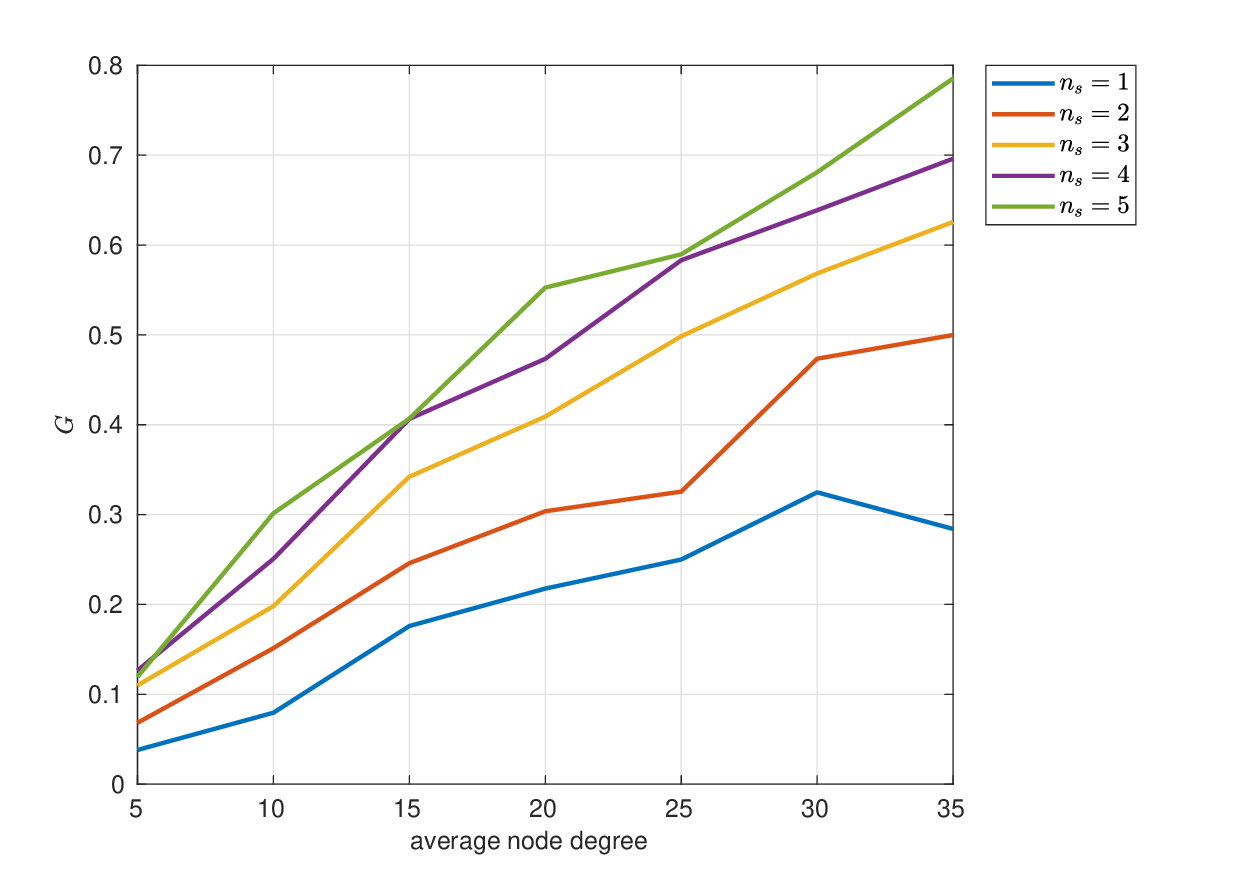}
  \caption{The cost of introducing the pro-rata rule.}\label{fig.gain}
\end{figure}

It seems natural that $G$ is growing as the graph is becoming more dense, since in this situation the pro-rata rule visibly reduces the number of free variables in the optimization problem. We have tested this conjecture using the random model described above.
The random graph contained $n=50$ nodes, whereas the average node degree $d=np$ varied from $0$ to $35$.
The number $n_s$ of nodes that receive the shock varies from 1 to 5.
The results were averaged over 50 runs. The resulting dependence between $G$ and $d$ is illustrated in Figure~\ref{fig.gain}.
We can see that the cost $G$ can be up to 79\%.
We can observe that, when there is a stronger shock, the relaxation of the pro-rata rule has a larger impact.

\begin{figure}
  \centering
  \includegraphics[width=0.62\columnwidth]{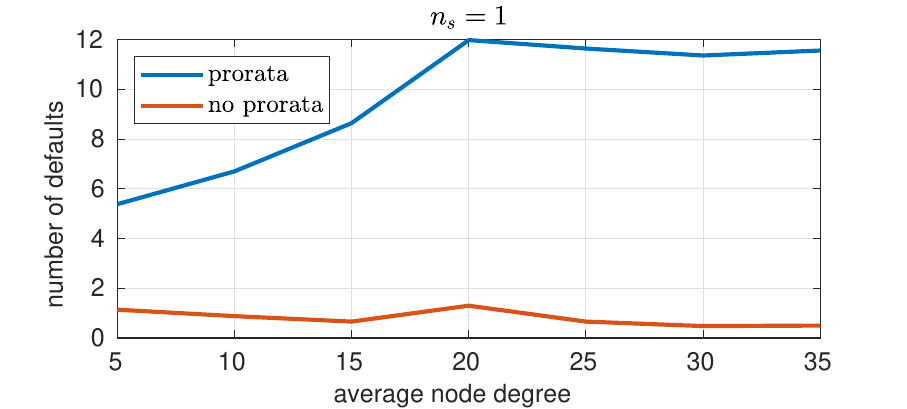}
  \includegraphics[width=0.62\columnwidth]{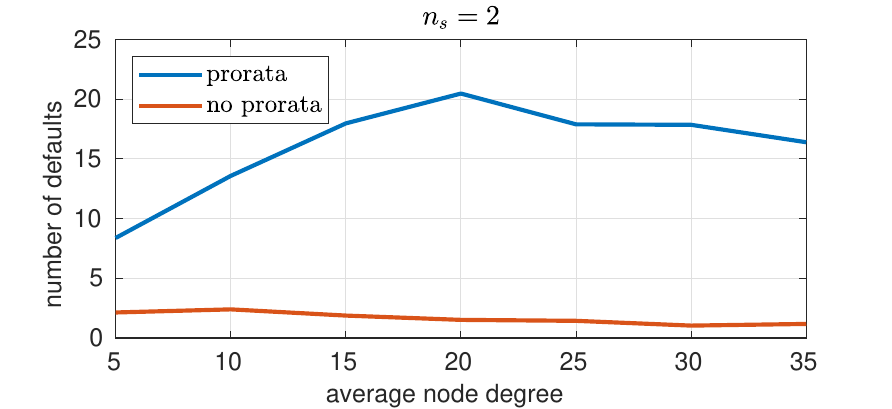}
  \includegraphics[width=0.62\columnwidth]{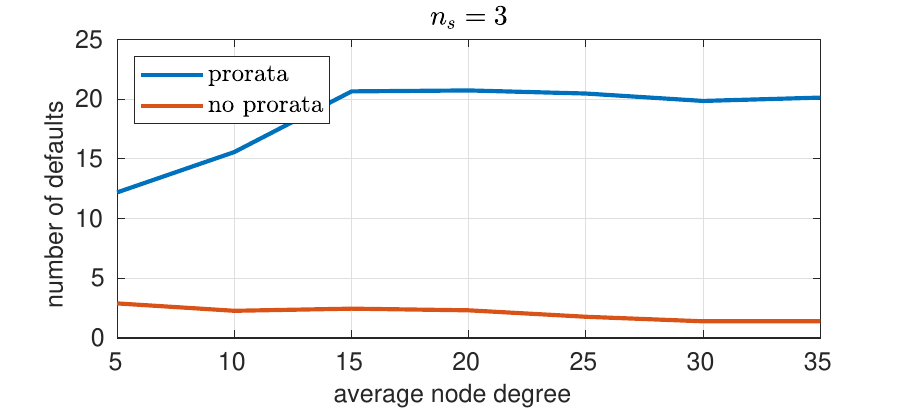}
  \includegraphics[width=0.62\columnwidth]{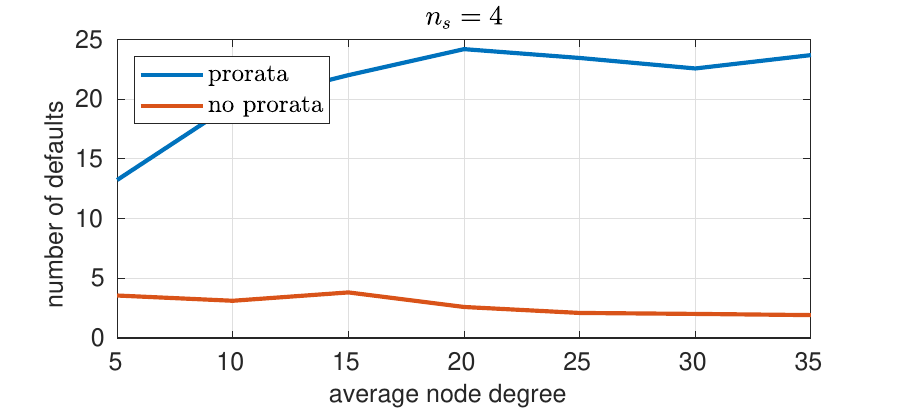}
  \includegraphics[width=0.62\columnwidth]{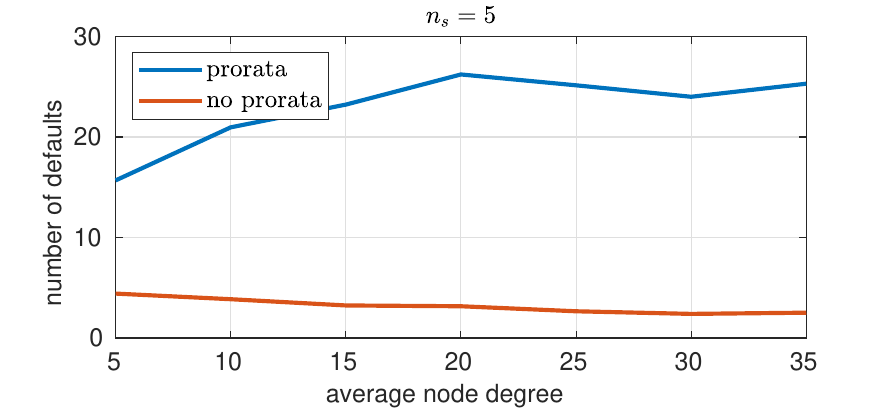}
  \caption{Numbers of defaulting nodes with or without the pro-rata rule.}
  \label{fig.defaults}
\end{figure}

{To evaluate the price of the pro-rata rule, we also introduce another metric that measures the number of defaulted nodes. This metric evaluates the dimension of the failure cascade caused by the initial shock. Figure \ref{fig.defaults} compares the number of defaulted nodes ($\phi\ap{out}_i<\bar\phi\ap{out}_i$) with or without the pro-rata rule. We  observe that relaxing the pro-rata rule significantly reduces the cascade failures and promotes isolation of the default contagion. }
  %\cred{[I would enrich the numerical examples as follows:
  %(1) make experiments with multiple injections of shocks, say $\alpha$ percent of the nodes receive a shock, etc.
  %(2) Evaluate a posteriori (i.e., after optimization) also other cost indices, such as the cardinality of the set of defaulted nodes.
  %(3) Maybe additional  examples with varying number of nodes, say 10, 20, 50, 100, 1000, and different nodes-to-degree ratios, to get an idea of how the performance and improvement varies with the network dimension...]}

\section{Conclusions}\label{sec:conclusion}
Based on the financial networks model of~\cite{Eisenberg2001}, we explored in this paper the concept of a clearing vector of payments, and
we developed new necessary and sufficient conditions for its uniqueness, together with a characterization of the set of all clearing vectors, see Theorem~\ref{thm.unique} and  Theorem~\ref{thm.unique1}.
Further,  we examined matrices of clearing payments that naturally arise if one relaxes the pro-rata rule.
Unique optimal clearing matrices can be computed efficiently via a two-stage  convex optimization approach.
Using numerical experiments with randomly generated synthetic networks, we showed that relaxation of the pro-rata rule allows to reduce
significantly  the aggregated  systemic loss and the total number of defaulted nodes, thus providing effective default isolation.

A few aspects remain critical for the proposed approach. First and foremost, an existing gap  from theory to practice is due to the fact  that
an aggregated unconstrained approach to the management of systemic risk should be accepted and contractualized ex-ante by the participant
players. Further, the inter-bank liability structure should be made transparent and available to a central regulatory authority who decides the clearing payments in the case of defaults. Alternatively, further research could be devoted to exploring the feasibility of dispensing from centralized knowledge of the whole liability matrix and of developing a system-level but distributed and decentralized optimal clearing payment algorithm whose iterations are based only on {\em local} exchange of information among neighboring nodes.

%\section{Acknowledgement}

\bibliographystyle{siamplain}
\bibliography{financial}

\begin{thebibliography}{10}

\bibitem{Acemoglu2015}
{\sc D.~Acemoglu, A.~Ozdaglar, and A.~Tahbaz-Salehi}, {\em Systemic risk and
  stability in financial networks}, American Economic Review, 105 (2015),
  pp.~564--608.

\bibitem{Allen2012}
{\sc F.~Allen, A.~Babus, and E.~Carletti}, {\em Asset commonality, debt
  maturity and systemic risk}, Journal of Financial Economics, 104 (2012),
  pp.~519--534.

\bibitem{Amini2023}
{\sc H.~Amini and Z.~Feinstein}, {\em Optimal network compression}, European
  Journal of Operational Research, 306 (2023), pp.~1439--1455.

\bibitem{Amini2016}
{\sc H.~Amini, D.~Filipovi{\'c}, and A.~Minca}, {\em To fully net or not to
  net: Adverse effects of partial multilateral netting}, Operations Research,
  64 (2016), pp.~1135--1142.

\bibitem{Barratt2020}
{\sc S.~Barratt and S.~Boyd}, {\em Multi-period liability clearing via convex
  optimal control}, Available at SSRN 3604618,  (2020).

\bibitem{Battiston2010}
{\sc S.~Battiston, J.~B. Glattfelder, D.~Garlaschelli, F.~Lillo, and
  G.~Caldarelli}, {\em The structure of financial networks}, in Network
  Science, Springer, 2010, pp.~131--163.

\bibitem{BlundellWignall2010}
{\sc A.~Blundell-Wignall and P.~Slovik}, {\em The {EU} stress test and
  sovereign debt exposures}, in {OECD} Working Papers on Finance, Insurance and
  Private Pensions, No. 4, OECD Financial Affairs Division, 2010,
  \url{https://www.oecd.org/finance/financial-markets/45820698.pdf}.

\bibitem{Calafiore2021a}
{\sc G.~Calafiore, G.~Fracastoro, and A.~Proskurnikov}, {\em On optimal
  clearing payments in financial networks}, in IEEE Conf. Decision and Control,
  2021, pp.~4804--4810.

\bibitem{Calafiore2022}
{\sc G.~Calafiore, G.~Fracastoro, and A.~Proskurnikov}, {\em Control of dynamic
  financial networks}, IEEE Control Systems Letters, 6 (2022), pp.~3206--3211.

\bibitem{Calafiore2023}
{\sc G.~C. Calafiore, G.~Fracastoro, and A.~Proskurnikov}, {\em Clearing
  payments in dynamic financial networks}, Automatica, 158 (2023), p.~111299.

\bibitem{Cifuentes2005}
{\sc R.~Cifuentes, G.~Ferrucci, and H.~S. Shin}, {\em Liquidity risk and
  contagion}, Journal of the European Economic Association, 3 (2005),
  pp.~556--566.

\bibitem{Csoka2018}
{\sc P.~Cs{\'o}ka and P.~Jean-Jacques~Herings}, {\em Decentralized clearing in
  financial networks}, Management Science, 64 (2018), pp.~4681--4699.

\bibitem{Eisenberg2001}
{\sc L.~Eisenberg and T.~H. Noe}, {\em Systemic risk in financial systems},
  Management Science, 47 (2001), pp.~236--249.

\bibitem{ElBitar2017}
{\sc K.~El~Bitar, Y.~M. Kabanov, and R.~Mokbel}, {\em On uniqueness of clearing
  vectors reducing the systemic risk}, Inform. Primen., 11 (2017),
  pp.~109--118.

\bibitem{Elliott2014}
{\sc M.~Elliott, B.~Golub, and M.~O. Jackson}, {\em Financial networks and
  contagion}, American Economic Review, 104 (2014), pp.~3115--53.

\bibitem{Elsinger2009}
{\sc H.~Elsinger et~al.}, {\em Financial networks, cross holdings, and limited
  liability}, Oesterreichische Nationalbank Austria, 2009.

\bibitem{Elsinger2006}
{\sc H.~Elsinger, A.~Lehar, and M.~Summer}, {\em Risk assessment for banking
  systems}, Management Science, 52 (2006), pp.~1301--1314.

\bibitem{Fischer2014}
{\sc T.~Fischer}, {\em No-arbitrage pricing under systemic risk: Accounting for
  cross-ownership}, Mathematical Finance: An International Journal of
  Mathematics, Statistics and Financial Economics, 24 (2014), pp.~97--124.

\bibitem{Gale2007}
{\sc D.~M. Gale and S.~Kariv}, {\em Financial networks}, American Economic
  Review, 97 (2007), pp.~99--103.

\bibitem{Gallice2019}
{\sc A.~Gallice}, {\em Bankruptcy problems with reference-dependent
  preferences}, International Journal of Game Theory, 48 (2019), pp.~311--336.

\bibitem{Glasserman2016}
{\sc P.~Glasserman and H.~P. Young}, {\em Contagion in financial networks},
  Journal of Economic Literature, 54 (2016), pp.~779--831.

\bibitem{Haldane2011}
{\sc A.~Haldane and R.~May}, {\em Systemic risk in banking ecosystems}, Nature,
  469 (2011), pp.~351--355.

\bibitem{Harary1965}
{\sc F.~Harary, R.~Norman, and D.~Cartwright}, {\em Structural models. An
  introduction to the theory of directed Graphs}, Wiley \& Sons, New York,
  London, Sydney, 1965.

\bibitem{Herings2021}
{\sc P.~J.-J. Herings and P.~Cs\'{o}ka}, {\em Uniqueness of clearing payment
  matrices in financial networks}, GSBE Research Memoranda, 014 (2021),
  pp.~1--29, \url{https://doi.org/10.26481/umagsb.2021014}.

\bibitem{Hurd2016}
{\sc T.~Hurd}, {\em Contagion! Systemic Risk in Financial Networks}, Springer
  Nature, 2016.

\bibitem{Jackson2021}
{\sc M.~O. Jackson and A.~Pernoud}, {\em Systemic risk in financial networks: a
  survey}, Annual Review of Economics, 13 (2021), pp.~171--202.

\bibitem{Kabanov2018}
{\sc Y.~M. Kabanov, R.~Mokbel, and K.~El~Bitar}, {\em Clearing in financial
  networks}, Theory of Probability \& Its Applications, 62 (2018),
  pp.~252--277.

\bibitem{Kaminski2000}
{\sc M.~M. Kaminski}, {\em 'hydraulic' rationing}, Mathematical Social
  Sciences, 40 (2000), pp.~131--155.

\bibitem{Kusnetsov2019}
{\sc M.~Kusnetsov and L.~A.~M. Veraart}, {\em Interbank clearing in financial
  networks with multiple maturities}, SIAM Journal on Financial Mathematics, 10
  (2019), pp.~37--67.

\bibitem{Massai2022}
{\sc L.~Massai, G.~Como, and F.~Fagnani}, {\em Equilibria and systemic risk in
  saturated networks}, Mathematics of Operation Research, 47 (2022),
  pp.~1707--2545.

\bibitem{Nier2007}
{\sc E.~Nier, J.~Yang, T.~Yorulmazer, and A.~Alentorn}, {\em Network models and
  financial stability}, Journal of Economic Dynamics and Control, 31 (2007),
  pp.~2033--2060.

\bibitem{Proskurnikov2020}
{\sc A.~Proskurnikov, G.~Calafiore, and M.~Cao}, {\em Recurrent averaging
  inequalities in multi-agent control and social dynamics modeling}, Annual
  Reviews in Control, 49 (2020), pp.~95--112.

\bibitem{Ren2016}
{\sc X.~Ren and L.~Jiang}, {\em Mathematical modeling and analysis of
  insolvency contagion in an interbank network}, Operations Research Letters,
  44 (2016), pp.~779--783.

\bibitem{Rogers2013}
{\sc L.~C. Rogers and L.~A. Veraart}, {\em Failure and rescue in an interbank
  network}, Management Science, 59 (2013), pp.~882--898.

\bibitem{Shin2008}
{\sc H.~S. Shin}, {\em Risk and liquidity in a system context}, Journal of
  Financial Intermediation, 17 (2008), pp.~315--329.

\bibitem{Suzuki2002}
{\sc T.~Suzuki}, {\em Valuing corporate debt: the effect of cross-holdings of
  stock and debt}, Journal of the Operations Research Society of Japan, 45
  (2002), pp.~123--144.

\bibitem{Thomson2013}
{\sc W.~Thomson}, {\em Game-theoretic analysis of bankruptcy and taxation
  problems: Recent advances}, International Game Theory Review, 15 (2013),
  p.~1340018.

\end{thebibliography}

 \section{Appendix}\label{sec:appendix}
%\vspace{.2cm}

\subsection{Technical preliminaries}

The proof of the main results of this paper are based on few technical propositions.
%
%please DON'T CORRECT 4.3a', prime is NOT a misprint here!!!
The next proposition follows, e.g., from~\cite[Corollary~4.3$a^\prime$]{Harary1965}
%please DON'T CORRECT 4.3a', prime is NOT a misprint here!!!
\begin{proposition}\label{prop.sink}
Each graph contains at least one sink strong component. Any node whose component is not a sink is a connected to one of the sink components by a path.
\end{proposition}

We will also use a technical lemma establishing necessary and sufficient conditions for the \emph{Schur stability} of substochastic matrices. The spectral radius of a square matrix $A$ is denoted by $\rho(A)$. We call a matrix $A$ \emph{Schur} stable if $\rho(A)<1$.
The Gershgorin disk theorem implies that $\rho(A)\leq 1$ for any substochastic matrix $A$. If $A$ is stochastic, then $\rho(A)=1$ and $A$ cannot be Schur stable.
\begin{lemma}\label{lem.substoch}
Let $A\in\Real{\calV\times\calV}$ be a substochastic matrix. Then three statements are equivalent:
\begin{enumerate}
\item Matrix $A$ is Schur stable: $\rho(A)<1$;
\item Submatrix $A^0=(a_{i,j})_{i,j\in \calV_0}$ is \emph{not} stochastic for every subset of indices $\calV_0\subseteq\calV$;
\item The set of nodes $\calV_d=\{i: \sum_{j}a_{ij}<1\}$ is non-empty and globally reachable in graph $\calG[A]$;
\item Each strong component of $\calG[A]$ being a sink contains at least one node from $\calV_d$.
\end{enumerate}
%In particular, an irreducible substochastic matrix is either Schur stable or stochastic.
\end{lemma}
 {\bf Proof of Lemma~\ref{lem.substoch}.}
Implication $3\Longrightarrow 1$ is a well-known fact, see e.g.~\cite[Lemma~7]{Proskurnikov2020}.
To prove 1$\Longrightarrow$2 notice that if submatrix $A^0=(a_{i,j})_{i,j\in \calV_0}$ is stochastic, then
 $0\leq\sum_{k\not\in\calV_0}a_{ik}\leq 1-\sum_{j\in\calV_0}a_{jk}=0\,\forall i\in\calV_0$, that is, $a_{ik}=0\,\forall i\in\calV_0,k\not\in\calV_0$. Hence, $A$ is
decomposable as follows
\beq\label{eq.decompose}
A=\begin{pmatrix}
A^0 & \mathbf{O}\\
* & *
\end{pmatrix},
\eeq
where $\mathbf{O}$ is a block of zeros and $*$ stand for some other submatrices. Hence, $1$ is an eigenvalue of $A$, in particular, $\rho(A)=1$ and $A$ is not Schur stable.

Implication 2$\Longrightarrow$4 is straightforward from the definition of a sink component. If $\calV_0$ is the set of nodes belonging to such a strong component, then matrix $A$ is decomposed as in~\eqref{eq.decompose}, where $A^0=(a_{i,j})_{i,j\in \calV_0}$. In view of statement 2, $A^0$ is not a stochastic matrix, so the sum in at least one of its rows is less than $1$, that is, $\calV_0\cap\calV_d\ne\emptyset$.

Implication 3$\Longrightarrow$4 is straightforward from Proposition~\ref{prop.sink} and the definition of a strongly connected component (whose every two nodes are mutually reachable).\qed

\subsection{Proofs of the main results}

{\bf Proof of Lemma~\ref{lem.optimum}.}
Since $\calD$ is compact (closed and bounded), the projection map $p\mapsto p_i$ reaches a maximal value $p_i^*=\max_{p\in\calD} p_i$. The vector $p^*\doteq (p_i^*)_{i\in\calV}$, by construction, dominates every vector from $\calD$. It remains to show that this vector belongs to $\calD$. By definition, $p_i^*\in [0,\bar p_i]$, $\forall i$.
Also, $c+A^{\top}p^*\geq c+A^{\top}p\geq p$, $\forall p\in\calD$. Hence, for each index $i\in\calV$ one has
\[
(c+A^{\top}p^*)_i\geq p_i\quad\forall p\in\calD.
\]
Taking the maximum over all $p$, one shows that
\[
(c+A^{\top}p^*)_i\geq p_i^*=\max_{p\in\calD} p_i,\;\forall i\in\calV,
\]
i.e., $c+A^{\top}p^*\geq p^*$, finishing the proof of statement 1.

Statement 2 is now straightforward from Definition~\ref{def.decrease}.

Statement 3 now follows from~\cite[Lemma~4]{Eisenberg2001}, stating that every minimizer in problem~\eqref{eq:clearing-opt-general} with a decreasing function $F$ is a clearing vector.
We give here the proof for the reader's convenience. Assume, by the purpose of contradiction, that~\eqref{eq:clearing-2} fails to hold, that is, $p_i^*<\min\left(\bar p_i,c_i+\sum_{k\ne i}a_{ki}p_k^*\right)$ for some $i$. Consider the vector $\hat p\doteq p^*+\delta\mathbf{e}_i$, where $\delta>0$ is sufficiently small and $\mathbf{e}_i$ is the coordinate vector whose $i$th coordinate is $1$ and others are null. By noticing that $\left(c + A\tran \hat p\right)_i>\hat p_i$ and
\[
\left(c + A\tran \hat p\right)_j\geq \left(c + A\tran p^*\right)_j\geq p_j^*=\hat p_j\quad\forall j\ne i,
\]
one shows that $\hat p\in\calD$, which contradicts  Statement~1.

To prove statement 4, denote the set of nodes of a sink component by $\calV_0$. Then, stochastic matrix $A$ decomposes as in~\eqref{eq.decompose}, and submatrix $A^0=(a_{i,j})_{i,j\in \calV_0}$ is also stochastic. Furthermore, this matrix is irreducible due to the definition of a strongly connected component. Thanks to the Perron-Frobenius theorem, $A^{\top}$ has a positive eigenvector $\pi^0$ such that $A\tran\pi^0=0$. Denoting
\[
\pi_i=\pi_i^0\,\forall i\in\calV_0,\quad \pi_i=0\,\forall i\not\in\calV_0,
\]
one obviously has $A\tran\pi=\pi$. Notice that vector $p^{\ve}=p^*+\ve\pi$ obeys the inequality $p^{\ve}\leq c+A\tran p^{\ve}$ for each $\ve>0$ and, furthermore, $p^{\ve}_i=p^*_i\leq\bar p_i\,\forall i\not\in\calV_0$. Since $p^*$ is the maximal element in $\calD$, $p^{\ve}\not\in\calD\,\forall\ve>0$, which is possible only when $p^*_i=\bar p_i$ for some $i\in\calV_0$.

Statement~5 follows from Lemma~\ref{lem.substoch} and the maximality of $p^*$. Let $p$ be a clearing vector enjoying the property from statement~4 and $I=\{i:p_i=\bar p_i\}$. Notice that, in view of Proposition~\ref{prop.sink}, the set $I$ is globally reachable in graph $\calG[A]$ (each node is connected to a sink strong component, and all nodes within this component are connected to some node from $I$). Denoting $I^c=\calV\setminus I$,  submatrix $\hat A=(a_{ij})_{i,j\in I^c}$, obviously, cannot contain stochastic submatrices, being thus Schur stable (Lemma~\ref{lem.substoch}).  The corresponding subvector $\hat p=(p_j)_{j\in I^c}$ obeys the equation
\beq\label{eq.aux}
\hat p=\hat A\tran p+\hat c,\quad \hat c_j\doteq c_j+\sum_{i\in I}a_{ij}\bar p_i\,\forall j\in I^c.
\eeq
Recalling that $p^*\geq p$, one has $p_i^*=\bar p_i\,\forall i\in I$; for this reason, subvector $\hat p^*=(p_j^*)_{j\in I^c}$ also obeys~\eqref{eq.aux}, which means that $\hat p=\hat p^*$ and, therefore, $p=p^*$.
\qed

\vspace{.2cm}
{\bf Proof of Lemma~\ref{lem.positivity}.}
We start with ``if'' part. If $c_i>0$ and $\bar p_i>0$, then polytope $\calD$, obviously, contains the vector $\ve\mathbf{e}_i$, where $\ve\in (0,\min(c_i,\bar p_i))$. Since $p^*=\max\calD$, we have $p_i^*>0$ whenever condition 1) holds. Notice also that if $p_j^*>0, p_i^*>0$ and $a_{ji}>0$, then~\eqref{eq:clearing-2} implies that $p_i^*>0$. Hence, condition 2) also entails that $p_i^*>0$. Assume now that condition 3) holds and let $\calI$ be the set of nodes of the strongly connected component containing $i$. Since $i$ is not a sink, this component is non-trivial (contains two or more nodes), and for each $j\in\calI$ we have
$\bar p_j>0$. Also, the submatrix $\tilde A=(a_{jk})_{j,k\in\calI}$ is stochastic and irreducible, hence, in view of the Perron-Frobenius theorem, there exist a \emph{strictly positive} eigenvector $v\in\mathbb{R}^{\calI}$ such that $A\tran v=v$. Rescaling, one may assume that $v_i\leq\bar p_i\,\forall i\in\calI$. Hence, $\calD$ contains the vector $p$, where
\[
p_j\doteq\begin{cases}
v_j, j\in\calI\\
0, j\not\in\calI
\end{cases}\quad\forall j\in\calV,
\]
in particular, $p_i^*\geq p_i>0$.

To prove the ``only if'' part,  notice first that $p_i^*>0$ entails that $\bar p_i>0$, so $i$ is not a sink node. Assume, on the contrary, that none of 1)-3) holds.
Consider again the strongly connected component $\calI\ni i$ and the corresponding irreducible  submatrix $\tilde A=(a_{jk})_{j,k\in\calI}$. Since condition 3) is violated, we have $\rho\doteq\rho(\tilde A)<1$.
Introducing the \emph{positive} Perron-Frobenius eigenvector $v\in\calI$ such
that $\tilde A\tran v=v$ and ~\eqref{eq:clearing-2} entails that
\beas
\lefteqn{
0<\sum_{j\in\calI} v_jp_j^*\leq \sum_{j\in\calI} v_j\sum_{k\in\calI}a_{kj}p_k^*+\sum_{j\in\calI} v_j\sum_{k\not\in\calI}a_{kj}p_k^*\leq}
\\
&& \leq\rho \sum_{k\in\calI} v_kp_k^*+\sum_{j\in\calI} v_j\sum_{k\not\in\calI}a_{kj}p_k^*.
\eeas
The latter inequality may hold only if $k\not\in\calI$ exists such that $a_{kj}p_k^*>0$. In other words, there exists another strongly connected component $\calI_1\ne\calI$ such that $p_j^*>0\,\forall j\in\calI_1$ and
$\calI_1$ is connected to $\calI$ by at least one arc (this implies that no arc can go from $\calI$ to $\calI_1$). Since condition 2) is violated for $i$, conditions 1 and 2 are also violated for every $i_1\in\calI_1$. Now we can repeat the argument, replacing $\calI$ by $\calI_1$ and find another strongly connected component $\calI_2\ne\calI_1$, such that at least one arc comes from $\calI_2$ to $\calI_1$ (in view of this $\calI_2\ne\calI$), $p_j^*>0\,\forall j\in\calI_2$ and none element of $\calI_2$ satisfies condition 1 or 2. Repeating this process, one could construct the infinite set of disjoint strongly connected components $\calI_0\doteq\calI,\calI_1,\calI_2,\ldots$ such that $\calI_s$ is connected to $\calI_{s-1}$, $p_j^*>0\,\forall j\in\calI_s$ yet conditions 1 and 2 fail to hold for $i\in\calI_s$. This leads to the contradiction (the set of nodes $\calV$ is finite). \qed

\vspace{.2cm}
{\bf Proof of Lemma~\ref{lem.suffic}.}
Consider an arbitrary clearing vector $p^0$ and let $J\doteq\{i:p_i^0<\bar p_i\}$. Our goal is to show that $p^0_i=p^*_i\,\forall i\in I_0'$.

\textbf{Step 1.} In view of~\eqref{eq:clearing-2}, the vector $p=p^0$ obeys the equations
\beq\label{eq.aux1}
\begin{gathered}
p_i=\bar p_i\quad\forall i\not\in J,\\
p_i=(c+A^{\top}p)_i=c_i+\sum_{k\not\in J}a_{ki}\bar p_k+\sum_{j\in J}a_{ji} p_j\;\forall i\in J.
\end{gathered}
\eeq

\textbf{Step 2.} We are going to show that if~\eqref{eq.aux1} has at least one solution, then $p=p^*$ should also be a solution to~\eqref{eq.aux1}. Indeed, every solution to system~\eqref{eq.aux1} is a global minimizer in the optimization problem~\eqref{eq:clearing-opt-general} with the objective function
\begin{gather}
F(p)\doteq\sum_{k\not\in J}(\bar p_k-p_k)+\sum_{j\in J}\left(c_j+\sum_{i\in\calV}a_{ij}p_i-p_j\right)\label{eq:phi}.
\end{gather}
If $p^0$ satisfies~\eqref{eq.aux1}, then $F(p^0)=0\leq F(p)\,\forall p\in\calD$. The function~\eqref{eq:phi} is non increasing, because it can be written as
\[
F(p)=\beta-b\tran p,\quad b_i\doteq 1-\sum_{j\in I}a_{ij}\geq 0,\quad\beta=\rm const,
\]
Lemma~\ref{lem.optimum} guarantees that $p^*$ is a global minimizer in~\eqref{eq:clearing-opt-general}, and hence $F(p^*)=0$. Recalling that $p^*\in\calD$, this is possible only if all addends in the right-hand side of~\eqref{eq:phi} vanish as $p=p^*$.

\textbf{Step 3.} We now show that the  equations~\eqref{eq.aux1} uniquely determine the subvector $(p_i)_{i\in I_0'}$. The first equation in~\eqref{eq.aux1} entails that $p_i=\bar p_i\,\forall i\in I_0'\setminus J$. By construction, the set $I_0'$ cannot be reached from any other node $j\not\in I_0'$ (otherwise, $I_0$ could also be reached from $j$), that is,
$a_{ji}=0$ whenever $i\in I_0'$ and $j\not\in I_0'$. Hence, for all solutions to~\eqref{eq.aux1} one has
\beq\label{eq.aux1a}
p_i=c_i+\sum_{k\not\in J}a_{ki}\bar p_k+\sum_{j\in I_0'\cap J}a_{ji} p_j\quad\forall i\in I_0'\cap J.
\eeq
We are now going to show that the matrix $\tilde A=(a_{ij})_{i,j\in I_0'\cap J}$ is Schur stable.
Assume, by purpose of contradiction , that $\tilde A$ is not Schur stable and thus (see Lemma~\ref{lem.substoch}) contains a \emph{stochastic} submatrix $(a_{ij})_{i,j\in K}$, where $K\subseteq I_0'\cap J$. Since $A$ is also stochastic, one has $a_{ij}=0\,\forall i\in K,j\in\calV\setminus K$. Denoting
\[\
\tilde z\doteq(z_i)_{i\in\calV},\quad
z_i\doteq\begin{cases}
1,\,i\in K,\\
0,\,i\not\in K,
\end{cases}
\]
one obviously has $Az\geq z$ or, equivalently, $z\tran A\tran\geq z\tran$.
The second equation in~\eqref{eq.aux1} entails that
\beas
\lefteqn{
z\tran p=\sum_{i\in K} z_ip_i=\sum_{i\in K} z_i\left(c_i+\sum_{j\in\calV}a_{ji}p_j\right)\overset{z_i=0\,\forall i\not\in K}{=\joinrel=\joinrel=\joinrel=} } \\
&& =\sum_{i\in\calV} z_i\left(c_i+\sum_{j\in\calV}a_{ji}p_j\right)\geq z\tran c+z\tran p,
\eeas
whence $z\tran c=0$ and thus $K\cap C^+=\emptyset$. On the other hand, $K\subseteq J\subseteq\calV\setminus S$, therefore, $K$ and $I_0$ are disjoint. By construction, graph $\calG$ contains no arcs connecting $K$ to nodes from $\calV\setminus K$, in particular, $I_0$ is not reachable from $K$. This contradicts to the definition of $I_0'\supseteq K$. The contradiction shows that $\tilde A$ is Schur stable, in particular,~\eqref{eq.aux1a} has a unique solution. Since both $p^0$ and $p^*$ satisfy~\eqref{eq.aux1} and~\eqref{eq.aux1a}, $p_i^0=p_i^*$ for all $i\in I_0'$.\qed

\vspace{.2cm}
{\bf Proof of Theorem~\ref{thm.unique}.}
The first statement is obvious, since the sets $I_0',I_1',\ldots$ are disjoint and the set of all nodes $\calV$ is finite.

To prove the second statement, we show it via induction on  $k=0,1,\ldots$  that for each clearing vector $p$ and each $k$ we have $p_i=p_i^*\,\forall i\in I_k'$. The induction base ($k=0$) is proved in Lemma~\ref{lem.suffic}.
Assume that the statement has been proved for $k=0,\ldots,q-1$ (where $q\geq 1$); we have to prove it for $k=q$.  By construction, banks from $\calV_q$ ($q\geq 1$) have no liability  to banks from $\calV\setminus\calV_q=I_0'\cup\ldots\cup I_{q-1}'$. Also, they neither have outside assets ($c_i=0\,\forall i\not\in I_0$) nor receive payments from banks from $I_0',\ldots,I_{q-2}'$. Only banks from $I_{q-1}'$ can pay liability to them. For each clearing vector $p$,~\eqref{eq:clearing-2} entails that
\beq
p_i=\min\left(\bar p_i,\sum_{s\in I_{q-1}'}a_{si}p_s+\sum_{k\in\calV_q}a_{ki}p_k\right)\quad\forall i\in\calV_q.\label{eq.aux3}
\eeq
In view of the induction hypothesis, the first sum is nothing else than $c_i^{(q)}$ for all clearing vectors $p$. Hence, for any clearing vector $p$ (of the original financial network) the subvector $(p_i)_{i\in\calV_q}$
serves a clearing vector in the \emph{reduced} financial network, determined by the set of nodes $\calV_q\ne\emptyset$, the respective submatrix $A^q=(a_{ij})_{i,j\in\calV_q}$ and the vector $c^{(q)}$, standing for the ``external'' assets.  This holds, in particular, for the dominant clearing vector $p=p^*$. Lemma~\ref{lem.suffic} guarantees that the elements $p_i,i\in I_q'$ are determined uniquely, that is, $p_i=p_i^*\,\forall i\in I_q'$, which proves the induction step $k=q$. This finishes the proof of the second statement of Theorem~\ref{thm.unique}. %from which the ``if'' part of the third statement is straightforward.

Furthermore, it is now obvious that each clearing vector (that is, a solution to~\eqref{eq:clearing-2}) has the structure~\eqref{eq.all-clearing-vec}, where $\xi=(p_i)_{i\in\calV_s}\geq 0$ is some vector. If $\calV_s=\emptyset$ (and the component $\xi$ is empty), the only clearing vector is $p=p^*$. Otherwise, $\xi$ obeys~\eqref{eq.subvector-constraint}. Indeed, the matrix $B$ is stochastic by construction, so that $\one\tran B=\one\tran$. In view of~\eqref{eq:clearing-2}, one has
\[
\xi\leq B\tran\xi.
\]
Multiplying the latter inequality by $\one\tran$, one shows that it can only hold if $\xi=B\tran\xi$. On the other had, $\xi_i\leq \bar p_i\,\forall i\in\calV_s$.
The inverse statement is also obvious: each vector~\eqref{eq.all-clearing-vec} satisfying~\eqref{eq.subvector-constraint} is a solution to~\eqref{eq:clearing-2}.
To finish the proof of the theorem, it remains to show that there are infinitely many subvectors $\xi$ obeying~\eqref{eq.subvector-constraint}. Notice that $B$ is a stochastic matrix, and  hence $\rho(B)=1$. Also, $\bar p_i>0\,\forall i\in\calV_s$ (the nodes cannot be sinks). The Perron-Frobenius theorem entails the existence of at least one non-negative eigenvector $\xi^0\in\Real{\calV_s}\setminus\{0\}$, which obeys~\eqref{eq.subvector-constraint}. Obviously, the set of all vectors satisfying~\eqref{eq.subvector-constraint} is a convex polytope, which contains a trivial vector $\xi=0$ and, thus, also the whole line segment $[0,\xi^0]$. We have proved that the set of clearing vectors in infinite.
\qed

\vspace{.2cm}
{\bf Proof of Theorem~\ref{thm.unique1}.}
Suppose that (i) holds and consider some non-trivial sink component whose set of nodes is $\calI$. Obviously, $\bar p_i>0\,\forall i\in\calI$.
The submatrix $\tilde A=(a_{ij})_{i,j\in\calI}$ is stochastic, let $v\in\mathbb{R}^{\calI}$ stand for the Perron-Frobenius eigenvector of $A\tran$ such that $A\tran v=v$. Assume that $\calI$ neither intersects $C^+$ nor is reachable from $C^+$. If $j\not\in\calI$ is connected to some node $i\in\calI$, then $c_j=0$ and $j$ cannot be reached from $C^+$. In view of Lemma~\ref{lem.positivity}, $p_j^*=0$ and, therefore, for all vectors from $\calD$ we have $p_j=0$. Therefore, every vector $p$ such that
\[
p=\begin{cases}
p_i^*,i\not\in\calI,\\
\ve v_i,i\in\calI
\end{cases}
\]
satisfies~\eqref{eq:clearing-2} when $\ve>0$ is small enough. We arrive at the contradiction with the assumption of the clearing vector's uniqueness. Implication (i)$\Longrightarrow$(ii) is proved.

Suppose now that (ii) holds yet the clearing vector is not unique. Consider the final set $\calV_s\ne\emptyset$ found by Algorithm~\ref{alg.1} ($s\geq 1$). By construction, no arc leads from $\calV_s$ to $I_0',\ldots,I_{s-1}'$, therefore, $\calV_s$ contains at least one sink component. This sink component cannot be trivial, because all sink nodes belong to $I_0\subseteq I_0'$. For the same reason,
this sink component contains no nodes from $C^+$. In view of (ii), a path from $C^+$ to $\calV_s$ exists. In view of Lemma~\ref{lem.positivity}, for each node $j$ on this path we have $p_j^*>0$.
Hence, there should exist an arc $j\to i$ connecting some $j\not\in\calV_s$ (with $p_j^*>0$) to some $i\in\calV_s$. This contradicts to the assumption that $\calV_s$ is a final set, because $c^{(s)}>0$.
Implication (ii)$\Longrightarrow$(i) is proved.
\qed

\vspace{.2cm}
{\bf Proof of Lemma~\ref{lem:opt_no_prorata}.}
Assume, by the purpose of contradiction,  that a local minimizer $P^*$ fails to satisfy~\eqref{eq:clearing-1}. Hence, an index $i$ exists such that
\[
p^*_i\doteq (P^*\one)_i=\sum_j p_{ij}^*<\min\left(\bar p_i,c_i+\sum_{k}p_{ki}^*\right).
\]
Since $p^*_i<\bar p_i$, there exists some $j$ such that $p^*_{ij}<\bar p_{ij}$. Then for small $\delta>0$
matrix $\hat P=P^*+\delta\mathbf{e}_i\mathbf{e}_j\tran$ (obtained from $P^*$ by the replacement
$p^*_{ij}\mapsto p^*_{ij}+\delta$, all other entries being invariant) belongs to $\calD_{n\times n}$.
Indeed, $\hat P\in [0,\bar P]$ and
\[
\hat p\doteq \hat P\one=p^*+\delta\mathbf{e}_i\leq c+(P^*)\tran\one\leq c+(\hat P)\tran\one.
\]
At the same time, $F(\hat P)<F(P^*)$. This contradicts to the assumption of local optimality. \qed

\vspace{.2cm}
{\bf Proof of Lemma~\ref{eq:lem_allsols}.}
The proof is immediate by observing that if $\tilde P$ is optimal for \eqref{eq_clearingopt_free_inflow1} then it must be that
$F_1(\tilde P) = F_1( P^*) $, where
\beas
F_1(\tilde P) &=&  \|\bar P - \tilde P\|_1 = \one\tran(\bar P - \tilde P)\one \quad\mbox{[since $\bar P \geq \tilde P$]} \\
F_1( P^*) &=&  \|\bar P - P^*\|_1 = \one\tran(\bar P -  P^*)\one \quad\mbox{[since $\bar P \geq  P^*$].}
\eeas
Therefore, $F_1(\tilde P) = F_1( P^*)$ implies that $\one\tran( \tilde P - P^*)\one = 0$. Defining $\Delta \doteq \tilde P - P^*$, we
have that $\tilde P = P^*+\Delta$, with $\Delta$ such that $\one\tran\Delta\one = 0$. Since $\tilde P$ must be feasible (i.e., it should satisfy
$0\leq  \tilde P \leq \bar P$, $c+\tilde P\tran \one -  \tilde P \one \geq 0$), we immediately obtain that $\tilde P\in \calS^* $.

Conversely, take any $\tilde P\in \calS^* $. Then, by construction, $\tilde P$ is feasible for \eqref{eq_clearingopt_free_inflow1} and
its objective is
\beas
F_1(\tilde P) &=&  \|\bar P - \tilde P\|_1 = \one\tran(\bar P - \tilde P)\one \quad\mbox{[since $\bar P \geq \tilde P$]} \\
&=&  \one\tran(\bar P -  P^* - \Delta)\one =  \one\tran(\bar P -  P^*) \one - \one\tran \Delta\one \\
&=& \one\tran(\bar P -  P^*) \one \quad \mbox{[since $\one\tran \Delta\one=0$ in the definition of $\calS^*$] }\\
&=& \|\bar P -  P^*\|_1 = F_1(P^*).
\eeas
Therefore, $\tilde P$ is optimal for \eqref{eq_clearingopt_free_inflow1}, which concludes the proof.
\qed

\vspace{.2cm}
{\bf Proof of Proposition~\ref{eq:prop_unique}.}
Problem  \eqref{prob:delta_unique}
is quadratic and strongly convex in its matrix variable $\Delta$. Further, this problem is feasible, since it admits at least the point
$\Delta = 0$, due to the fact that $P^*$ is optimal, hence feasible, for \eqref{eq_clearingopt_free_inflow1}.
Therefore, $\Delta^*$ exists and it is unique.
Consider then  $\tilde P^* \doteq P^* + \Delta^*$. Since by construction $\bar P \geq \tilde P^* \geq 0$ and
$c+ \tilde P^{*\top}\one -  \tilde P^{*}\one \geq 0$, we have that $\tilde P^*$ is feasible for problem \eqref{eq_clearingopt_free_inflow1}.
Furthermore, its objective value is
\beas
F_1(\tilde P^*) &=&  \|\bar P - \tilde P^*\|_1 = \one\tran(\bar P - \tilde P^*)\one \quad\mbox{[since $\bar P \geq \tilde P^*$]} \\
&=&  \one\tran(\bar P -  P^* - \Delta^*)\one =  \one\tran(\bar P -  P^*) \one - \one\tran \Delta^*\one \\
&=& \one\tran(\bar P -  P^*) \one \quad \mbox{[since $\one\tran \Delta^*\one=0$ due to the constraint in \eqref{prob:delta_unique}] }\\
&=& \|\bar P -  P^*\|_1 = F_1^*.
\eeas
Therefore, $\tilde P^*$ is optimal for \eqref{eq_clearingopt_free_inflow1}, which proves point (b) in Proposition~\ref{eq:prop_unique}.
 From Lemma~\ref{lem:opt_no_prorata}
we also immediately  conclude that $\tilde P^*$ is a clearing matrix, which proves point (a).
For point (c), consider that from Lemma~\ref{eq:lem_allsols} it holds that $P^* +\Delta$ in problem \eqref{prob:delta_unique}
spans over all optimal solutions to problem  \eqref{eq_clearingopt_free_inflow1}. Therefore,  problem \eqref{prob:delta_unique}
provides the unique minimum Euclidean norm solution among all the possible optimal solutions to
 \eqref{eq_clearingopt_free_inflow1}.
\qed

\end{document}